\documentclass[11pt,a4paper]{article}
\usepackage[margin=27mm]{geometry}
\usepackage[numbers,sort&compress]{natbib}

\usepackage[T1]{fontenc}
\usepackage{lmodern}
\usepackage{amsmath,amssymb,amsthm}
\usepackage{graphicx}
\usepackage{booktabs}
\usepackage{microtype}
\usepackage{xcolor}
\definecolor{linkblue}{RGB}{0,70,150}
\usepackage[colorlinks=true,citecolor=linkblue,linkcolor=linkblue,urlcolor=linkblue]{hyperref}
\hypersetup{pdftitle={Error Estimates and Convergence Analysis for Learning Mechanism-AI Coupled Dynamical Systems},
 pdfauthor={Pengfei Song; Yanni Xiao; Yuan Lou; Shuangshuang Yin}}
\newcommand{\U}{\mathcal{U}}
\newcommand{\X}{\mathcal{X}}

\theoremstyle{plain}
\newtheorem{theorem}{Theorem}
\newtheorem{lemma}{Lemma}
\newtheorem{proposition}{Proposition}

\newtheorem{assumption}{Assumption}[section]
\theoremstyle{definition}
\newtheorem{definition}{Definition}
\theoremstyle{remark}
\newtheorem{remark}{Remark}

\title{Error Estimates and Convergence Analysis for Learning Mechanism--AI Coupled Dynamical Systems}
\author{\normalsize Pengfei Song$^{1}$\quad Yanni Xiao$^{1}$\quad Yuan Lou$^{2}$\quad Shuangshuang Yin$^{1,*}$}
\date{}
\begin{document}
\raggedbottom
\maketitle
\begin{center}
\small
$^{1}$Department of Applied Mathematics, School of Mathematics and Statistics,\\
Xi'an Jiaotong University, Xi'an 710049, China\\[3pt]
$^{2}$School of Mathematical Sciences, CMA-Shanghai,\\
Shanghai Jiao Tong University, Shanghai 200240, China\\[5pt]
\href{mailto:song921012@xjtu.edu.cn}{song921012@xjtu.edu.cn};
\href{mailto:yxiao@mail.xjtu.edu.cn}{yxiao@mail.xjtu.edu.cn};
\href{mailto:yuanlou@sjtu.edu.cn}{yuanlou@sjtu.edu.cn}\\
$^{*}$Corresponding author:
\href{mailto:yinss2014@stu.xjtu.edu.cn}{yinss2014@stu.xjtu.edu.cn}
\end{center}
\begin{abstract}
Mechanism--AI coupled dynamical systems retain known mechanisms while learning unknown dynamics from data. We study how neural approximation, numerical integration, finite observations, and optimization jointly affect the accuracy of this learning process. For time-dependent learned terms in a fixed compact admissible class, we derive a four-component objective error bound within an optimal control formulation. Explicit estimates are obtained for constrained ReLU approximation and Euler integration, with deterministic quadrature and random sampling treated separately. Under joint refinement, we establish $\Gamma$-convergence of the approximate objectives and show that approximate minimizers approach the continuous optimal set when the global optimization gaps vanish. We also derive conditional stability estimates for state-dependent learned terms. Numerical experiments with a susceptible--infected--recovered (SIR) model illustrate the analysis and show how solver resolution affects continuous-trajectory accuracy even when training losses are similar.
\end{abstract}
\noindent\textbf{Keywords:} Optimal control; Universal differential equations; Error estimates; Learning theory; Mechanism-AI coupled dynamic model.

\medskip
\noindent\textbf{Mathematics Subject Classification:} 37M99; 65L20; 68T07; 62M45.

\section{Introduction}
Combining mechanistic models with data-driven terms provides a way to study biological systems whose governing processes are only partly known. This approach retains established physical or biological principles while inferring unresolved dynamics from observations. Among these hybrid approaches, the formulation of \emph{universal differential equations} (UDEs) \cite{rackauckas2020universal,yin2023optimal}
\begin{equation}\label{eq0}
	\frac{\mathrm{d}x}{\mathrm{d}t} = f_{\eta}(x, t, \text{UniversalApprox}_{\phi}(x, t))
\end{equation}
encodes known mechanisms in differential equations. Unknown components, such as control terms, external forces, or closure relations, are represented by trainable function approximators, typically neural networks. Here $x$ is the state, $\eta$ denotes fixed mechanistic parameters, and $\phi$ denotes the learned parameters. Our main results concern the time-dependent case $u(t)$ in ordinary differential equations (ODEs). Section~\ref{sec2} gives a separate stability estimate for the state-dependent case in \eqref{eq0} and identifies the additional approximation and numerical assumptions needed for a full extension.
In \cite{yin2023optimal}, the learning problem was formulated as an optimal control task, where unknown control functions are parameterized by neural networks and optimized through differentiable numerical solvers. This formulation raises a theoretical question: under what conditions do the resulting approximate solutions approach an optimum of the continuous problem? Convergence to an optimal loss must be distinguished from recovery of the underlying dynamics, which requires identifiability from the available observations. To address the convergence question, we analyze approximation, numerical, generalization, and optimization errors.

Optimal control and convergence analyses of neural differential equations have substantial precedents. Schwartz and Polak \cite{schwartz1996consistent} developed consistent approximations of constrained optimal control problems based on Runge--Kutta integration. Their variational approach connects discrete optimization problems with their continuous counterparts. Hager \cite{hager2000runge} analyzed discretization orders for optimal controls and adjoint variables, which need not equal the forward integration order. For NeurODE training, Bonnet et al.\ \cite{bonnet2023measure} developed a mean-field optimal control analysis, including finite-sample stability and quantitative generalization estimates. Aghili and Mula \cite{aghili2024adaptive} studied the interaction between sampling and time-discretization errors and used these estimates to design adaptive procedures.

The separation of approximation, estimation, and optimization errors is also established in learning theory \cite{bottou2008tradeoffs,shen2022fixed}. Jentzen and Welti \cite{jentzen2023overall} give an overall error analysis for a specified neural empirical-risk minimization scheme trained by stochastic gradient descent with random initialization. In the setting of dynamical systems, Kang and Gong \cite{kang2022feedforward} study approximation and complexity for compositional functions and numerical algorithms, including differential equations and optimal control. Dong et al.\ \cite{dong2022optimization} analyze optimization problems whose differential-equation constraints contain learned components. Scagliotti \cite{scagliotti2023deep} uses optimal control and $\Gamma$-convergence to study flow-based approximation of diffeomorphisms and generalization.

Our analysis concerns constrained, time-dependent terms coupled to a prescribed mechanistic ODE. We place admissible neural approximation, uniform solver error, and statistical estimation of the numerical trajectory loss in the same control space. Explicit bounds for constrained ReLU interpolation and Euler integration verify these requirements for a concrete class. Under joint refinement and vanishing global objective gaps, the finite-domain estimates yield convergence of minimum values and ensure that approximate minimizers approach the continuous optimal set. The SIR example connects this analysis with coefficient identification for positive initial susceptible and infected populations and with the effect of solver resolution on fitted trajectories.

We analyze approximations of an optimal control problem in \textit{Lagrange} form \cite{betts2010practical,lenhart2007optimal}:
\begin{equation}\label{e1}
	\left\{
	\begin{aligned}
		& \min_{u(\cdot),x(\cdot)} J = \int_{0}^{T}g(x,u,t)dt                              \\
		& s.t. \quad \frac{\mathrm{d}x}{\mathrm{d}t} = f(x,u,t), x(0) = x_{0}, t\in [0,T], \\
		& \quad \quad u_{i}(t) \in [a_{i}, b_{i}], 1\leq i \leq m, t\in [0,T], \\
& \quad \quad \operatorname{Lip}_{\infty}(u)\le L_0,
	\end{aligned}
	\right.
\end{equation}
Here $T>0$, $a_i\le b_i$ are finite, $u:[0,T]\to\mathbb R^m$, $f:\mathbb R^n\times\mathbb R^m\times[0,T]\to\mathbb R^n$, and $g:\mathbb R^n\times\mathbb R^m\times[0,T]\to\mathbb R$. The initial state $x_0$ and mechanistic parameters are fixed.
Terminal costs can be treated by extending the objective; we restrict attention to the integral cost in \eqref{e1}. In addition to the range constraints, we impose a fixed regularity bound $0\le L_0<\infty$. With $A=\prod_{i=1}^m[a_i,b_i]$, the admissible controls are
\begin{equation}\label{eq:admissible_controls}
 \mathcal U=\{u\in C([0,T];A):\operatorname{Lip}_{\infty}(u)\le L_0\},
 \qquad d(u,v)=\|u-v\|_\infty.
\end{equation}
Here $\|u\|_\infty=\sup_{t\in[0,T]}|u(t)|_\infty$ and $\operatorname{Lip}_{\infty}$ uses the same vector norm. Thus \eqref{e1} is studied on this specified subclass of range-constrained controls. The common Lipschitz bound is maintained as the network complexity increases. The class is compact in the stated metric by the Arzel\`a--Ascoli theorem. State trajectories lie in a common compact convex invariant set $K\subset\mathbb R^n$, and we use $\|x\|_{\mathcal X}=\sup_{t\in[0,T]}|x(t)|_1$. Section~\ref{sec2} states the assumptions on $f$ and $g$ that guarantee well-posedness and continuity in these norms.
For any fixed $u \in \mathcal{U}$, the ordinary differential equation $\frac{\mathrm{d}x}{\mathrm{d}t} = f(x,u,t), x(0) = x_{0}, t\in [0,T]$ admits a unique solution $S(u)$, and we call $S:\mathcal{U}\rightarrow \mathcal{X}$ the control-to-state map. Under the assumptions stated below, the continuous objective attains its minimum on $\mathcal U$. An optimal control, which need not be unique, is denoted by
\begin{equation}
	u^{*} \in \arg\min_{u\in \mathcal{U}}J(S(u),u).
\end{equation}
To approximate the optimal control problem \eqref{e1}, we parameterize an admissible control $u(t)$ by a neural network
\begin{equation}\label{e2}
	u(t) = \mathrm{NeuralNetwork}_{\theta}(t),
\end{equation}
thereby restricting $u$ to a surrogate class $\mathcal{U}_{N}\subset\mathcal U$ with a prescribed architecture and complexity $N\in\mathbb{Z}^{+}$. Examples include shallow networks with rectified linear unit (ReLU) activation, deep networks of fixed width, and feedforward networks represented by directed acyclic graphs (DAGs) \cite{kang2022feedforward}. Depending on the architecture, $N$ may count the width, depth, or total number of units. The feasible parameter set $\Theta_N$ is nonempty and compact, its realization map $\theta\mapsto u_\theta$ is continuous in the uniform norm, and its outputs satisfy \eqref{eq:admissible_controls}. These conditions ensure that $\mathcal U_N$ is a nonempty compact class, rather than a linear subspace. A constraint-preserving ReLU construction is given in Section~\ref{sec2}. The resulting parameter optimization problem is:
\begin{equation}\label{e3}
	\left\{
	\begin{aligned}
		& \min_{\theta\in\Theta_N} J = \int_{0}^{T}g(x,\mathrm{NeuralNetwork}_{\theta}(t),t)dt                                         \\
		& s.t. \quad \frac{\mathrm{d}x}{\mathrm{d}t} = f(x,\mathrm{NeuralNetwork}_{\theta}(t),t), x(0) = x_{0}, t\in [0,T], \\
	\end{aligned}
	\right.
\end{equation}
and a minimizer of this restricted problem is denoted by
\begin{equation}\label{e4}
	u_N^{*}\in\arg\min_{u\in \mathcal{U}_N}J(S(u),u).
\end{equation}
This restriction introduces \emph{approximation error}: the distance from an optimal control to its best approximation in the chosen admissible neural class. The approximation must preserve both range and regularity constraints. Section~\ref{sec2} verifies this property for a nodal ReLU construction and derives its error bound.
Evaluating the objective in \eqref{e3} requires solving the state equation \cite{wanner1993solving,wanner1996solving}. For fixed $u\in \mathcal{U}_N$, denote the numerical solution of $\frac{\mathrm{d}x}{\mathrm{d}t} = f(x,u,t), x(0) = x_{0}, t\in [0,T]$ by
$$x_{h}=S_{h}(u).$$
Define
\begin{equation}\label{e5}
	J^{h}(S_{h}(u),u)=
	\left\{
	\begin{aligned}
		& J(S_{h}(u),u), u \in \mathcal{U}_N,        \\
		& +\infty, u \in \mathcal{U} \setminus \mathcal{U}_N,
	\end{aligned}
	\right.
\end{equation}
where $h$ is the maximum integration step and $S_h(u)$ denotes a continuous reconstruction of the numerical states. We assume that these reconstructed states remain in $K$. A minimizer of the numerical problem is denoted by
\begin{equation}\label{e6}
	u_{N,h}^{*}\in\arg\min_{u\in \mathcal{U}_N}J^{h}(S_{h}(u),u).
\end{equation}
Replacing $S(u)$ by $S_h(u)$ introduces \emph{numerical error}. A uniform state-error estimate and the Lipschitz continuity of the cost integrand control its effect on the objective. The required stability and consistency estimates depend on the solver and the regularity of the admissible controls.
To evaluate \eqref{e6} during optimization, the integral in $J^{h}(S_h(u),u)$ can be approximated by numerical quadrature. For nodes $t_i\in[0,T]$ and nonnegative weights $\xi_i$ with $\sum_{i=1}^{\delta}\xi_i=T$, define
\begin{equation}\label{e7}
	J^{h,\delta}(S_{h}(u), u)=
	\left\{
	\begin{aligned}
		& \sum_{i=1}^{\delta} \xi_{i} g(S_{h}(u)(t_{i}), u(t_{i}), t_{i}), u \in \mathcal{U}_N, \\
		& +\infty, u \in \mathcal{U} \setminus \mathcal{U}_N,
	\end{aligned}
	\right.
\end{equation}
For a data-driven loss, we instead observe independent pairs $Z_i=(t_i,Y_i)$ drawn from the same distribution, with $t_i$ uniform on $[0,T]$. Let $\ell(x,u,Y)$ be a bounded measurable observation loss, uniformly Lipschitz in $(x,u)$. The population integrand is $g(x,u,t)=\mathbb E[\ell(x,u,Y)\mid t]$, and the empirical objective is
\begin{equation}\label{e8}
	J^{h,\delta}(S_{h}(u), u)=
	\left\{
	\begin{aligned}
		& \frac{T}{\delta}\sum_{i=1}^{\delta}\ell(S_{h}(u)(t_{i}), u(t_{i}), Y_i), u \in \mathcal{U}_N, \\
		& +\infty, u \in \mathcal{U} \setminus \mathcal{U}_N,
	\end{aligned}
	\right.
\end{equation}
The factor $T$ makes this empirical objective consistent with the unnormalized integral in \eqref{e1}. Quadrature and random sampling are two alternative constructions of the finite objective. In either case, finding a restricted optimal control gives the nonlinear programming problem
\begin{equation}\label{e9}
	u_{N,h,\delta}^{*} \in \arg\min_{u\in \mathcal{U}_N}J^{h,\delta}(S_{h}(u), u).
\end{equation}
The statistical construction introduces \emph{generalization error}, the discrepancy between empirical and population objectives. Rademacher complexity \cite{shalev2014understanding} and uniform concentration control this discrepancy for the loss class induced by the numerical trajectories. The deterministic construction \eqref{e7} introduces \emph{quadrature error}, governed by the rule and the regularity of the integrand. We treat these two forms of loss evaluation separately.
Training for \eqref{e9} can use automatic differentiation to evaluate derivatives \cite{baydin2018automatic,ma2021comparison}. Parameter updates can then use first-order methods such as Adam (\cite{kingma2014adam}, \cite[Chapter 5]{kochenderfer2019algorithms}) or quasi-Newton methods such as the limited-memory BFGS method (L-BFGS) (\cite{jorge2006numerical}, \cite[Chapter 6]{kochenderfer2019algorithms}), with the admissibility constraints enforced. Let the feasible control produced at iteration $k$ be
\begin{equation}\label{e10}
	u_{N,h,\delta}^{k}\in\mathcal U_N,\qquad k=0,1,\ldots.
\end{equation}
For the reduced objectives, write
\begin{equation}\label{e12}
	\hat{J}(u)=J(S(u),u), \hat{J}^{h}(u)=J^{h}(S_{h}(u),u), \hat{J}^{h,\delta}(u)=J^{h,\delta}(S_{h}(u),u).
\end{equation}
The fourth component, \emph{optimization error}, is the gap between the objective value at the terminated iterate and the global minimum of the finite objective. Neural training is generally nonconvex, so convergence to a stationary point alone does not control this gap. We retain it explicitly and use an algorithmic rate only under assumptions that give a global objective estimate.

\section{Main Theoretical Results}
We first state a total error bound under the regularity assumptions of Section~\ref{sec2}. The bound separates approximation, numerical, generalization, and optimization errors. For deterministic evaluation of the integral, quadrature error replaces generalization error. We use the state and control norms introduced with \eqref{eq:admissible_controls}. The factor $T$ in the cost constants comes from integrating bounds in the uniform trajectory norm.

\begin{theorem}[Total Error Bound]\label{thm:error_bound}
Let $u^*$ be an exact optimal control for \eqref{e1} on $\mathcal U$, and let $u^k_{N,h,\delta}\in\mathcal U_N$ be a feasible control obtained after $k$ optimization steps. Assume the compactness, stability, and cost regularity conditions of Section~\ref{sec2}. Define
\[
 \mathrm{Aprox}_{u^*}(N)=\min_{v\in\mathcal U_N}\|v-u^*\|_\infty,
 \qquad
 \mathrm{Num}(h)=\sup_{v\in\mathcal U}\|S_h(v)-S(v)\|_{\mathcal X}.
\]
For the random observation model \eqref{e8}, fix $N$, $h$, and $\delta$ independently of the observations, and suppose that the bounded-loss and parameter-covering assumptions of Proposition~\ref{prop:loss_covering} hold. For any fixed $r>0$ and $\alpha\in(0,1)$, with probability at least $1-\alpha$,
\[
 0\le \hat J(u^k_{N,h,\delta})-\hat J(u^*)
 \le \mathcal E_{\mathrm{approx}}+\mathcal E_{\mathrm{num}}
      +\mathcal E_{\mathrm{gen}}+\mathcal E_{\mathrm{opt}},
\]
where the components defined in \eqref{e14} satisfy
\[
\begin{aligned}
 \mathcal E_{\mathrm{approx}}
   &\le T(L_{gx}L_S+L_{gu})\,\mathrm{Aprox}_{u^*}(N),\\
 \mathcal E_{\mathrm{num}}&\le 2TL_{gx}\,\mathrm{Num}(h),\\
 \mathcal E_{\mathrm{gen}}&\le 2T\zeta_{\delta,N,h}(r,\alpha),\\
 \mathcal E_{\mathrm{opt}}&=\rho_{N,h,\delta}(k)
   :=\hat J^{h,\delta}(u^k_{N,h,\delta})
       -\min_{v\in\mathcal U_N}\hat J^{h,\delta}(v).
\end{aligned}
\]
Here $\zeta_{\delta,N,h}$ is the explicit deviation bound in \eqref{eq:loss_covering}. The same event controls every feasible output in $\mathcal U_N$, including one selected by training on these observations.
For the deterministic quadrature objective \eqref{e7}, the inequality holds deterministically with the third term denoted by $\mathcal E_{\mathrm{quad}}$ and bounded by
\[
 \mathcal E_{\mathrm{quad}}\le 2\mathrm{Quad}(\delta),\qquad
 \sup_{v\in\mathcal U_N}
 |\hat J^{h,\delta}(v)-\hat J^h(v)|\le\mathrm{Quad}(\delta).
\]
The quadrature bound is taken uniformly over the solver meshes under consideration.
\end{theorem}

The four terms have the following meanings:
\begin{itemize}
 \item \textbf{Approximation error $\mathcal E_{\mathrm{approx}}$}: the error due to representing the control by a finite-capacity neural network. The quantity $\mathrm{Aprox}_{u^*}(N)$ measures the distance between an optimal control $u^*$ and its best admissible network approximation. For the constrained piecewise-linear ReLU class constructed below, it is at most $L_0T/(2N)$ and tends to zero as $N\to\infty$.
 \item \textbf{Numerical error $\mathcal E_{\mathrm{num}}$}: the error introduced by discretizing the differential equation. The function $\mathrm{Num}(h)$ is a uniform forward-state error, rather than a relative error or an error in the optimized control. Euler integration gives an explicit $O(h)$ bound under the common Lipschitz assumptions. Higher orders require the corresponding uniform regularity and solver conditions.
 \item \textbf{Generalization error $\mathcal E_{\mathrm{gen}}$}: the statistical error arising from a finite sample of size $\delta$. The relevant hypothesis class consists of observation losses evaluated on numerical ODE trajectories. Rademacher complexity bounds its uniform deviation, while the parameter-covering estimate below gives an explicit confidence bound. In the deterministic setting, the third term is instead the quadrature error of the specified rule.
 \item \textbf{Optimization error $\mathcal E_{\mathrm{opt}}$}: the global objective gap of the terminated iterate for the finite problem. If a particular algorithm and objective satisfy $\rho_{N,h,\delta}(k)\le C_{\mathrm{opt}}\mathrm{Opt}(k)$, that estimate can be substituted. Increasing $k$ alone does not establish that this gap tends to zero for a nonconvex neural objective.
\end{itemize}
Consequently, the total objective error tends to zero along a sequence for which the approximation, numerical, loss-evaluation, and global optimization errors all vanish. The error bound in Theorem~\ref{thm:error_bound} identifies what must be controlled when network capacity, integration accuracy, sample size, and optimization effort are varied together. It does not prescribe a universal convergence rate for a training algorithm.

We next establish the \textit{variational convergence} of the discretized optimal control problems to the continuous problem. The theory of $\Gamma$-convergence \cite{braides2002gamma,brezis1993introduction} then gives convergence of minimum values and characterizes the limit points of minimizing controls.

Choose a sequence $N_j\to\infty$, $h_j\to0$, and $\delta_j\to\infty$, and write $\epsilon_j=\max\{h_j,1/\delta_j,1/N_j\}$. The functional $\hat J_{\epsilon_j}:\mathcal U\to\mathbb R\cup\{+\infty\}$ equals the finite objective on $\mathcal U_{N_j}$ and $+\infty$ elsewhere. Assume compactness of $\mathcal U$ in the uniform metric, continuity of $\hat J$, consistency of the solver, admissible neural approximation, and vanishing uniform error on the finite domains, as detailed in Assumptions~\ref{ass:compact}--\ref{ass:uniform}.
\begin{theorem}[$\Gamma$-Convergence and Convergence of Minimizers]\label{thm:gamma_conv}
Under these assumptions, the sequence $\{\hat J_{\epsilon_j}\}$ is equi-coercive and $\Gamma$-converges to $\hat J$ in the uniform topology. Consequently,
\begin{enumerate}
 \item The minimum values converge:
 \[
 \min_{u\in\mathcal U}\hat J(u)
   =\lim_{j\to\infty}\min_{u\in\mathcal U}\hat J_{\epsilon_j}(u).
 \]
 \item Every sequence $\{u^*_{\epsilon_j}\}$ of discrete minimizers is precompact, and every limit point is a minimizer of $\hat J$.
\end{enumerate}
The same limit-point conclusion holds for feasible approximate minimizers whose global objective gaps tend to zero. For random observation objectives, the conclusions hold on any sample path on which the uniform error assumption holds. Under the bounded, uniformly Lipschitz observation-loss assumptions of Section~\ref{sec3}, using prefixes of one infinite iid observation sequence supplies such an event of probability one.
\end{theorem}

Theorem~\ref{thm:gamma_conv} establishes consistency of the specified discretization and admissible approximation. It complements the quantitative estimates of Theorem~\ref{thm:error_bound}. Compactness and uniform approximation on the finite feasible domains are sufficient for the proof; $\Gamma$-convergence records the liminf and recovery-sequence properties when those domains change. If the continuous minimizer is unique, the whole sequence of approximate minimizers converges to it. Without uniqueness or an identifiability condition, the theorem gives convergence to the optimal set and does not identify a particular physical coefficient. The subsequent sections prove these results and illustrate their implications with the SIR model.

\begin{remark}[The scope of this work.]
The results concern a fixed class of mechanism-AI coupled ODE control problems. The control constraints, neural realization, observation model, and solver estimates are chosen so that their compatibility can be verified. The common regularity bound supplies the compactness used in the proof; allowing this bound to grow with $N$ requires a separate analysis.

Optimization problems constrained by partial differential equations (PDEs) \cite{hinze2008optimization} raise related questions for neural-network-based approximations. Extending the present argument would require suitable state-space compactness, boundary conditions, uniform solver estimates, and admissible neural approximation. These properties are not consequences of the ODE analysis and are outside the results proved here.
\end{remark}

\section{Proof of Total Error Bound}\label{sec2}
We prove the error decomposition and total error bound in Theorem~\ref{thm:error_bound} under the stated Lipschitz assumptions.
The uniform error estimates will also be used in the subsequent $\Gamma$-convergence analysis. In that argument, the discrepancy between the exact objective $\hat{J}$ and its numerical counterpart $\hat{J}^{h,\delta}$ is estimated on the finite-valued domain $\mathcal{U}_N$.
We collect the controls introduced above:
\begin{equation}\label{e11}
	\begin{aligned}
		& u^{*} \in \arg\min_{u\in \mathcal{U}}J(S(u),u),\\
		& u_N^{*}\in\arg\min_{u\in \mathcal{U}_N}J(S(u),u),\\
		& u_{*}\in\arg\min_{u\in \mathcal{U}_N}\|u-u^{*}\|_{\mathcal{U}},\\
		& u_{N,h}^{*}\in\arg\min_{u\in \mathcal{U}_N}J^{h}(S_{h}(u),u),\\
		& u_{N,h,\delta}^{*}\in\arg\min_{u\in \mathcal{U}_N}J^{h,\delta}(S_{h}(u),u),\\
		& u_{N,h,\delta}^{k}\in\mathcal{U}_N.
	\end{aligned}
\end{equation}
All minimizers in \eqref{e11} are global. Their existence follows from the compactness and continuity conditions below. The last line denotes a computed admissible iterate; convergence to a global minimizer is not presumed.

We use the admissible control class
\[
\begin{gathered}
\mathcal{U}=\{u\in C([0,T];A): |u(t)-u(s)|_\infty\leq L_0|t-s|
\text{ for all }s,t\in[0,T]\},\\
A=\prod_{i=1}^m[a_i,b_i].
\end{gathered}
\]
Here $L_0\geq0$ is fixed, and the control and state norms are
\[
\begin{gathered}
\|u\|_{\mathcal{U}}=\sup_{t\in[0,T]}|u(t)|_\infty,\qquad
\|x\|_{\mathcal{X}}=\sup_{t\in[0,T]}|x(t)|_1,\\
\mathcal{X}=C([0,T];\mathbb{R}^n).
\end{gathered}
\]
The notation $\|\cdot\|_{\mathcal{U}}$ denotes the ambient uniform norm and also applies to differences of controls. The constrained set $\mathcal{U}$ is compact in this norm by the Arzel\`a--Ascoli theorem: its elements are uniformly bounded and equicontinuous, and its defining constraints are preserved under uniform limits. It is a subset of a Banach space, not itself a vector space.

We make the following assumptions in this section:
\begin{enumerate}
	\item The exact and numerical control-to-state maps take values in a common compact convex set $K\subset\mathbb{R}^n$ containing $x_0$. They are Lipschitz continuous: there exist constants $L_S,L_{S,h}\geq0$ such that for all $u,v\in\mathcal{U}$,
	\[
	\|S(u)-S(v)\|_{\mathcal{X}}\leq L_S\|u-v\|_{\mathcal{U}},\qquad
	\|S_h(u)-S_h(v)\|_{\mathcal{X}}\leq L_{S,h}\|u-v\|_{\mathcal{U}}.
	\]
	Sufficient conditions for the exact map are given below; a uniform numerical estimate is proved for Euler's method in the Numerical Error subsection. The condition that numerical trajectories remain in $K$ must be checked for the solver and step sizes used.
	\item In the deterministic setting, $g$ is continuous on $K\times A\times[0,T]$ and Lipschitz continuous in both state and control variables: for all $x_1,x_2\in K$, $u_1,u_2\in A$, and $t\in[0,T]$,
	\[
	|g(x_1,u_1,t)-g(x_2,u_2,t)|
	\leq L_{gx}|x_1-x_2|_1+L_{gu}|u_1-u_2|_\infty.
	\]
	For the statistical risk formulation, continuity in $t$ may be replaced by measurability and an integrable bound uniform in $(x,u)$, while retaining this Lipschitz estimate. The bounded-loss assumption below supplies such a bound. These conditions ensure that the integral objective is finite.
\end{enumerate}
The neural classes $\mathcal{U}_N$ are assumed nonempty and compact in $\mathcal{U}$, as verified for a concrete construction in the Approximation error subsection. Then the continuous objectives attain the minima in \eqref{e11}; for $\hat{J}^{h,\delta}$ this applies to each fixed quadrature rule or realized finite sample with a continuous loss.

\begin{lemma}\label{lemma1}
	\noindent (i) $\hat{J}$ is Lipschitz continuous on $\mathcal{U}$ with constant
	$C_J=T(L_{gx}L_S+L_{gu})$.

	\noindent (ii) The finite-valued restriction of $\hat{J}^h$ to $\mathcal{U}_N$ is Lipschitz continuous with constant $T(L_{gx}L_{S,h}+L_{gu})$.
\end{lemma}
To prove Lemma \ref{lemma1}, recall the definitions:
\[
\hat{J}(u)=J(S(u),u)=\int_0^T g(S(u)(t),u(t),t)\,dt,\qquad u\in\mathcal{U},
\]
\[
\hat{J}^h(u)=J^h(S_h(u),u)=\int_0^T g(S_h(u)(t),u(t),t)\,dt,\qquad u\in\mathcal{U}_N,
\]
where $S$ and $S_h$ are the exact and numerical control-to-state maps, respectively. The extension of $\hat{J}^h$ by $+\infty$ outside $\mathcal{U}_N$ is not part of the Lipschitz assertion.

\begin{proof}
(i) For any $u,v\in\mathcal{U}$, we have
	\begin{align*}
		|\hat{J}(u)-\hat{J}(v)|
		&=\left|\int_0^T[g(S(u)(t),u(t),t)-g(S(v)(t),v(t),t)]\,dt\right|\\
		&\leq\int_0^T|g(S(u)(t),u(t),t)-g(S(v)(t),v(t),t)|\,dt.
	\end{align*}
Using the Lipschitz continuity of $g$ and integrating over $[0,T]$ yields
	\begin{align*}
		|\hat{J}(u)-\hat{J}(v)|
		&\leq L_{gx}\int_0^T|S(u)(t)-S(v)(t)|_1\,dt\\
		&\quad+L_{gu}\int_0^T|u(t)-v(t)|_\infty\,dt\\
		&\leq T L_{gx}\|S(u)-S(v)\|_{\mathcal{X}}
		+T L_{gu}\|u-v\|_{\mathcal{U}}\\
		&\leq T(L_{gx}L_S+L_{gu})\|u-v\|_{\mathcal{U}}.
	\end{align*}
The factor $T$ follows from the uniform norms and the unnormalized integral objective.

(ii) The proof for $\hat{J}^h$ is identical, replacing $S$ with $S_h$ and $L_S$ with $L_{S,h}$. For any $u,v\in\mathcal{U}_N$,
	\begin{align*}
		|\hat{J}^h(u)-\hat{J}^h(v)|
		&\leq L_{gx}\int_0^T|S_h(u)(t)-S_h(v)(t)|_1\,dt\\
		&\quad+L_{gu}\int_0^T|u(t)-v(t)|_\infty\,dt\\
		&\leq T L_{gx}\|S_h(u)-S_h(v)\|_{\mathcal{X}}
		+T L_{gu}\|u-v\|_{\mathcal{U}}\\
		&\leq T(L_{gx}L_{S,h}+L_{gu})\|u-v\|_{\mathcal{U}}.
	\end{align*}
This completes the proof.
\end{proof}

\begin{remark}
The constants $L_{gx}$ and $L_{gu}$ can be derived from gradient bounds of $g$. For example, if $g$ is $C^1$ on a neighborhood of $K\times A\times[0,T]$, the dual norms corresponding to the chosen state and control norms give
	\[
	L_{gx}=\sup_{K\times A\times[0,T]}|\nabla_x g|_\infty,\qquad
	L_{gu}=\sup_{K\times A\times[0,T]}|\nabla_u g|_1.
	\]
These suprema are finite by continuity and compactness.
\end{remark}

\begin{remark}
The Lipschitz constants $L_S$ and $L_{S,h}$ depend on the dynamical system and numerical solver, respectively. To ensure the existence and regularity of the exact state map, suppose that $f(x,u,t)$ is continuous and locally Lipschitz in $x$ on a neighborhood of $K$, uniformly in $(u,t)\in A\times[0,T]$, and that $K$ is positively invariant for every $u\in\mathcal{U}$. On $K\times A\times[0,T]$, assume
	\[
	\begin{aligned}
	|f(x_1,u_1,t)-f(x_2,u_2,t)|_1
	&\leq L_{fx}|x_1-x_2|_1+L_{fu}|u_1-u_2|_\infty,\\
	|f(x,u,t)|_1&\leq M_f.
	\end{aligned}
	\]
Under these conditions, the state equation
	\[
	\dot{x}(t)=f(x(t),u(t),t),\qquad x(0)=x_0,
	\]
admits a unique solution on $[0,T]$ for each $u\in\mathcal{U}$. The common compact invariant set prevents finite-time escape, and each trajectory is $M_f$-Lipschitz in time. Subtracting the integral state equations and applying Gronwall's inequality gives
	\[
	\begin{aligned}
	|S(u)(t)-S(v)(t)|_1
	&\leq L_{fu}\int_0^t e^{L_{fx}(t-s)}|u(s)-v(s)|_\infty\,ds\\
	&\leq L_{fu}\Phi(L_{fx},t)\|u-v\|_{\mathcal{U}},
	\end{aligned}
	\]
where
	\[
	\Phi(c,t)=
	\begin{cases}(e^{ct}-1)/c,&c>0,\\t,&c=0.\end{cases}
	\]
Thus we may take $L_S=L_{fu}\Phi(L_{fx},T)$. For a discretized map $S_h$, a corresponding Lipschitz estimate requires numerical stability; consistency alone does not establish it.
\end{remark}

\begin{remark}[State-dependent learned terms]
For a feedback term $u=\nu(t,x)$, let $\nu:D\to A$ be $L_0$-Lipschitz, where
$D=[0,T]\times K$ has metric
$d_\infty((t,x),(s,y))=\max\{|t-s|,|x-y|_1\}$ and the output norm is $|\cdot|_\infty$.
Assume that the closed-loop equations are well posed and remain in $K$ for all such feedbacks.
For $\|\nu-\widetilde{\nu}\|_{\infty,D}
=\sup_{(t,x)\in D}|\nu(t,x)-\widetilde{\nu}(t,x)|_\infty$,
the same argument gives a state-map constant
	\[
	L_S^{\mathrm{fb}}=L_{fu}\Phi(L_{fx}+L_{fu}L_0,T).
	\]
Indeed, the difference of feedback values along two trajectories is at most
$\|\nu-\widetilde{\nu}\|_{\infty,D}+L_0|x-\widetilde{x}|_1$.
Consequently, the cost $J(x_\nu,\nu(\cdot,x_\nu))$ is Lipschitz in the feedback norm with constant
	\[
	C_J^{\mathrm{fb}}
	=T\big[(L_{gx}+L_{gu}L_0)L_S^{\mathrm{fb}}+L_{gu}\big].
	\]
This gives the change in the stability constants; admissible approximation classes and numerical consistency must also be established when applying the error decomposition to feedback networks.
\end{remark}

\subsection{Error Decomposition}
The error decomposition follows the four stages of the learning problem: choosing a hypothesis class, solving the dynamics numerically, evaluating the loss by quadrature or sampling, and optimizing the parameters. Adding and subtracting intermediate objective values and using their minimizing properties separates the errors estimated below. The same comparison will also support the convergence argument.

All the terms below are finite: $\hat{J}$ is evaluated on $\mathcal{U}$, and each numerical objective is evaluated only on $\mathcal{U}_N$. We obtain the following comparison:
\begin{equation}\label{e13}
	\begin{aligned}
		&0\leq\hat{J}(u_{N,h,\delta}^{k})-\hat{J}(u^{*})\\
		&\leq|\hat{J}(u_{N,h,\delta}^{k})-\hat{J}^{h}(u_{N,h,\delta}^{k})|
		+|\hat{J}^{h}(u_{N,h,\delta}^{k})-\hat{J}^{h,\delta}(u_{N,h,\delta}^{k})|\\
		&\quad+|\hat{J}^{h,\delta}(u_{N,h,\delta}^{k})-\hat{J}^{h,\delta}(u_{N,h,\delta}^{*})|\\
		&\quad+\hat{J}^{h,\delta}(u_{N,h,\delta}^{*})-\hat{J}^{h,\delta}(u_{N,h}^{*})
		\qquad(\leq0)\\
		&\quad+|\hat{J}^{h,\delta}(u_{N,h}^{*})-\hat{J}^{h}(u_{N,h}^{*})|\\
		&\quad+\hat{J}^{h}(u_{N,h}^{*})-\hat{J}^{h}(u_{N}^{*})
		\qquad(\leq0)\\
		&\quad+|\hat{J}^{h}(u_{N}^{*})-\hat{J}(u_{N}^{*})|\\
		&\quad+\hat{J}(u_{N}^{*})-\hat{J}(u_{*})
		\qquad(\leq0)\\
		&\quad+|\hat{J}(u_{*})-\hat{J}(u^{*})|\\
		&\leq\underbrace{|\hat{J}(u_{N,h,\delta}^{k})-\hat{J}^{h}(u_{N,h,\delta}^{k})|}_{\text{Numerical error}}\\
		&\quad+\underbrace{|\hat{J}^{h}(u_{N,h,\delta}^{k})-\hat{J}^{h,\delta}(u_{N,h,\delta}^{k})|}_{\text{Quadrature or generalization error}}\\
		&\quad+\underbrace{|\hat{J}^{h,\delta}(u_{N,h,\delta}^{k})-\hat{J}^{h,\delta}(u_{N,h,\delta}^{*})|}_{\text{Optimization error}}\\
		&\quad+\underbrace{|\hat{J}^{h,\delta}(u_{N,h}^{*})-\hat{J}^{h}(u_{N,h}^{*})|}_{\text{Quadrature or generalization error}}\\
		&\quad+\underbrace{|\hat{J}^{h}(u_{N}^{*})-\hat{J}(u_{N}^{*})|}_{\text{Numerical error}}\\
		&\quad+\underbrace{|\hat{J}(u_{*})-\hat{J}(u^{*})|}_{\text{Approximation error}}.
	\end{aligned}
\end{equation}
The three nonpositive differences follow from the definitions of
$u_{N,h,\delta}^{*}$, $u_{N,h}^{*}$, and $u_N^{*}$, respectively.
Define the four error components by
\begin{equation}\label{e14}
	\begin{aligned}
		\mathcal{E}_{\text{approx}}
		&=|\hat{J}(u_{*})-\hat{J}(u^{*})|,\\
		\mathcal{E}_{\text{gen}}
		&=|\hat{J}^{h}(u_{N,h,\delta}^{k})-\hat{J}^{h,\delta}(u_{N,h,\delta}^{k})|\\
		&\quad+|\hat{J}^{h,\delta}(u_{N,h}^{*})-\hat{J}^{h}(u_{N,h}^{*})|,\\
		\mathcal{E}_{\text{num}}
		&=|\hat{J}(u_{N,h,\delta}^{k})-\hat{J}^{h}(u_{N,h,\delta}^{k})|
		+|\hat{J}^{h}(u_{N}^{*})-\hat{J}(u_{N}^{*})|,\\
		\mathcal{E}_{\text{opt}}
		&=|\hat{J}^{h,\delta}(u_{N,h,\delta}^{k})-\hat{J}^{h,\delta}(u_{N,h,\delta}^{*})|
		=\rho_{N,h,\delta}(k).
	\end{aligned}
\end{equation}
The last equality uses the global objective gap $\rho_{N,h,\delta}(k)$ defined in Theorem~\ref{thm:error_bound}. For random observations, $\mathcal{E}_{\text{gen}}$ denotes the statistical loss-evaluation error. For deterministic quadrature, the same two comparison terms are denoted by $\mathcal E_{\mathrm{quad}}$. Their respective bounds are derived separately below.

\subsection{Approximation error}
Let $\mathcal{U}_{N}$ be the set of admissible networks with a specified architecture (e.g., shallow neural networks with ReLU activation, deep neural networks with fixed width, feedforward neural networks with a fixed finite number of neurons, or directed acyclic graphs \cite{kang2022feedforward}) and complexity $N\in\mathbb{Z}^{+}$ (e.g., width, depth, or the total number of units).
We require a nonempty compact parameter set $\Theta_N\subset\mathbb{R}^{P_N}$ and a continuous realization map $\theta\mapsto u_\theta$ from $\Theta_N$ into $\mathcal{U}$ in the uniform norm. Then $\mathcal{U}_N=\{u_\theta:\theta\in\Theta_N\}$ is compact. The amplitude and Lipschitz constraints are part of this admissible parameterization.
For any function $u\in\mathcal{U}$, denote the degree of approximation by
\begin{equation}\label{e15}
	\mathrm{Aprox}_{u}(N):=\mathrm{dist}(u,\mathcal{U}_{N})
	=\inf_{v\in\mathcal{U}_{N}}\|u-v\|_{\mathcal{U}}.
\end{equation}
The approximation property needed here is that, for every $u\in\mathcal{U}$,
\begin{equation}\label{e16}
	\mathrm{Aprox}_{u}(N)\rightarrow0\qquad(N\rightarrow\infty).
\end{equation}
Nestedness of the classes is not required. Compactness of $\mathcal{U}_N$, rather than finite dimensionality alone, ensures that
$\arg\min_{v\in\mathcal{U}_N}\|u-v\|_{\mathcal{U}}$ is nonempty, so $u_*$ is well defined. These network classes are generally nonlinear subsets of the function space. Thus Lemma \ref{lemma1} gives
\begin{equation}\label{e22}
	\mathcal{E}_{\text{approx}}
	=|\hat{J}(u_*)-\hat{J}(u^*)|
	\leq C_J\,\mathrm{Aprox}_{u^*}(N),
	\qquad C_J=T(L_{gx}L_S+L_{gu}).
\end{equation}

A concrete admissible construction is obtained by linear interpolation on $N$ equal time intervals. Set $\Delta_N=T/N$ and $t_j=j\Delta_N$. Choose nodal parameters
\[
\Theta_N=\{(c_0,\ldots,c_N)\in A^{N+1}:
|c_{j+1}-c_j|_\infty\leq L_0\Delta_N,\ 0\leq j<N\},
\]
and let $d_j=(c_{j+1}-c_j)/\Delta_N$. The interpolant has the exact ReLU representation
\[
u_\theta(t)=c_0+d_0t+\sum_{j=1}^{N-1}(d_j-d_{j-1})(t-t_j)_+,
\qquad (r)_+=\max\{r,0\}.
\]
The nodal set is compact and nonempty. Convexity of $A$ and the slope constraints give $u_\theta\in\mathcal{U}$, while
$\|u_\theta-u_{\widetilde{\theta}}\|_{\mathcal{U}}\leq
\max_j|c_j-\widetilde{c}_j|_\infty$ proves continuity of the realization map.
For $u\in\mathcal{U}$, choosing $c_j=u(t_j)$ produces an admissible interpolant $I_Nu$. On a cell, write $t=(1-\lambda)t_j+\lambda t_{j+1}$. Then
\[
|u(t)-I_Nu(t)|_\infty
\leq(1-\lambda)L_0\lambda\Delta_N
+\lambda L_0(1-\lambda)\Delta_N
\leq\frac{L_0\Delta_N}{2}.
\]
Consequently,
\[
\mathrm{Aprox}_{u}(N)\leq\|u-I_Nu\|_{\mathcal{U}}
\leq\frac{L_0T}{2N}.
\]
If $u\in W^{2,\infty}(0,T;\mathbb{R}^m)$, the linear-interpolation remainder gives the stronger bound
\[
\mathrm{Aprox}_{u}(N)\leq
\frac{T^2}{8N^2}\|u''\|_{L^\infty},
\]
where the derivative norm is the essential supremum of the componentwise maximum norm. Both estimates preserve the amplitude and slope constraints. Here $N$ counts time intervals; for fixed $m$, the number of network units and nodal parameters grows linearly in $N$. The equal-mesh classes need not be nested, but they satisfy \eqref{e16} along the full sequence.

For other network architectures, $\mathrm{Aprox}_{u^*}(N)$ depends on the target class, the norm, and how network complexity is counted:
\begin{itemize}
	\item Shallow ReLU networks: For the scalar time variable considered here, the preceding construction gives $O(N^{-1})$ for Lipschitz targets and $O(N^{-2})$ for targets in $W^{2,\infty}$. Related shallow-network bounds, such as Barron's \cite{barron2002universal}, concern sigmoidal networks, a Fourier-moment condition, and an integrated squared-error criterion. Their $L^2$ error rates cannot be substituted directly for the uniform ReLU approximation error in \eqref{e15}. Sobolev approximation results for ReLU networks \cite{yarotsky2017error} also specify the depth and complexity allowed; they do not give an arbitrary-order rate for a fixed shallow architecture.

	\item Deep ReLU networks: For deep networks, the relation between width, depth, parameter selection, and the target's modulus of continuity affects the attainable rate \cite{yarotsky2018optimal}. Montanelli and Yang \cite{montanelli2020error} obtain bounds using Kolmogorov--Arnold superposition structure and assumptions on the approximation of the outer functions. These hypotheses must be checked before using the resulting rate for a particular control class; no exponential rate for general piecewise smooth controls is assumed here.

	\item Graph-based architectures: For control functions defined on graph domains (e.g., in multi-agent systems), the stability result of Gama et al.\ \cite{gama2020stability} concerns relative perturbations of the graph. It provides a stability analysis, rather than a general $O(N^{-1/2})$ approximation theorem. Applying the present argument to a graph-based class requires a suitable admissible domain, norm, and a separate verification of \eqref{e16}.
\end{itemize}
In our analysis, we do not require a specific functional form of the approximation rate. We require \eqref{e16} within the constrained class $\mathcal{U}$; unrestricted universal approximation alone does not verify that requirement. The explicit admissible rates above can be substituted into Theorem \ref{thm:error_bound} to obtain a concrete approximation contribution.

\subsection{Numerical Error}
We quantify numerical error through an a priori estimate for the control-to-state map. The decomposition applies to any numerical method that satisfies this estimate.
For $x=S(u)$ and $x_h=S_h(u)$, define the uniform absolute state error
\begin{equation}\label{e17}
	\mathrm{Num}(h):=\sup_{u\in\mathcal{U}}
	\|S(u)-S_h(u)\|_{\mathcal{X}},\qquad
	\|x-x_h\|_{\mathcal{X}}\leq\mathrm{Num}(h),
\end{equation}
and require $\mathrm{Num}(h)\to0$ as $h\to0$.
The supremum ensures that consistency is uniform as the network class changes. Using the Lipschitz bound for $g$, we obtain
\begin{equation}\label{e18}
	\begin{aligned}
		\mathcal{E}_{\text{num}}
		&=|\hat{J}(u_{N,h,\delta}^{k})-\hat{J}^{h}(u_{N,h,\delta}^{k})|
		+|\hat{J}^{h}(u_N^*)-\hat{J}(u_N^*)|\\
		&\leq2T L_{gx}\,\mathrm{Num}(h).
	\end{aligned}
\end{equation}

For a sufficiently smooth initial value problem, the global error of a stable single-step method is of order $O(h^p)$, where $p$ is the method's order of convergence \cite{wanner1993solving,wanner1996solving}. The one-step local error is usually one order higher. To use such a rate in \eqref{e17}, the regularity and stability constants, and the accuracy of the continuous reconstruction, must hold uniformly over the admissible controls.
\begin{itemize}
	\item Explicit Euler Method: This method has $\mathrm{Num}(h)=O(h)$ under the common Lipschitz bounds and invariant-set condition specified below. The estimate applies to the continuous piecewise linear reconstruction.
	\item Runge--Kutta Methods (RK4): The classical fourth-order method has $O(h^4)$ state error when the required derivatives of the reduced vector field $f(x,u(t),t)$ are uniformly bounded on each integration interval and the continuous reconstruction has the corresponding accuracy. For piecewise smooth controls, the intervals must resolve the breakpoints, or the crossing steps need a separate error analysis. The amplitude and Lipschitz constraints alone do not imply fourth-order convergence.
	\item Implicit Methods (e.g., Backward Euler, Trapezoidal Rule): These methods are useful for stiff systems. Backward Euler is first-order, while the Trapezoidal Rule (also called Crank--Nicolson in the corresponding PDE setting) is second-order and A-stable. Their use in \eqref{e17} additionally requires well-defined implicit steps, common state bounds, and the uniform regularity appropriate to the stated order.
\end{itemize}

For completeness, a uniform Euler estimate can be proved directly with the regularity used here. In addition to the state and control Lipschitz bounds, assume
\[
|f(x,u,t)-f(x,u,s)|_1\leq L_{ft}|t-s|.
\]
Let $0=s_0<\cdots<s_M=T$, $h_j=s_{j+1}-s_j$, and $h=\max_jh_j$. Use the same mesh for all controls in the following comparison. Euler's method is
\[
x_{j+1}^h=x_j^h+h_j f(x_j^h,u(s_j),s_j),\qquad x_0^h=x_0,
\]
and $S_h(u)$ is its linear interpolant. Assume that the Euler iterates remain in $K$ for every $u\in\mathcal{U}$ and every mesh under consideration. Convexity of $K$ then also places the interpolation in $K$.

Set $D_f=L_{fx}M_f+L_{fu}L_0+L_{ft}$. The exact trajectory is $M_f$-Lipschitz, so $t\mapsto f(S(u)(t),u(t),t)$ is $D_f$-Lipschitz. The local residual is therefore at most $D_fh_j^2/2$. The nodal errors $e_j=|S(u)(s_j)-x_j^h|_1$ satisfy
\[
e_{j+1}\leq(1+L_{fx}h_j)e_j+\frac{D_fh_j^2}{2},\qquad e_0=0.
\]
Since $h_j^2\leq h h_j$ and
\[
\sum_{j=0}^{M-1}h_j\prod_{i=j+1}^{M-1}(1+L_{fx}h_i)
\leq\Phi(L_{fx},T),
\]
with the value $T$ when $L_{fx}=0$, discrete Gronwall gives
\[
\max_je_j\leq(D_fh/2)\Phi(L_{fx},T).
\]
The exact trajectory has Lipschitz derivative with constant $D_f$, so interpolating its exact nodal values adds at most $D_fh^2/8$. Hence
\[
\mathrm{Num}(h)\leq E_h:=
D_f\left[\frac{h}{2}\Phi(L_{fx},T)+\frac{h^2}{8}\right].
\]
The same discrete comparison, applied to two controls, gives
\[
\|S_h(u)-S_h(v)\|_{\mathcal{X}}\leq
L_{fu}\Phi(L_{fx},T)\|u-v\|_{\mathcal{U}},\qquad
\operatorname{Lip}_t(S_h(u))\leq M_f.
\]
The time-Lipschitz bound follows from the slopes of the linear reconstruction. Thus $L_{S,h}$ can be chosen as $L_S$, uniformly in the Euler mesh. The condition on numerical invariance remains part of this estimate and must be verified in each application.

An order estimate for the state equation is distinct from an order estimate for discretized optimal controls. The latter also involves the optimality system and its stability; Runge--Kutta order conditions in that setting are analyzed in \cite{hager2000runge}.

\subsection{Quadrature  error}
For $\hat{J}^{h,\delta}$ defined by numerical quadrature, write
\[
\varphi_u^h(t)=g(S_h(u)(t),u(t),t),\qquad
\hat{J}^{h,\delta}(u)=\sum_{i=1}^{\delta}\xi_i\varphi_u^h(t_i),
\qquad u\in\mathcal{U}_N.
\]
We use positive weights with $\sum_i\xi_i=T$.
The Bramble--Hilbert lemma \cite{bramble1970estimation} gives a higher-order estimate when the quadrature error annihilates polynomials of the required degree. For an integer $r\geq0$, a fixed composite rule exact through degree $r$, with a fixed number of nodes per cell and a quasiuniform refinement, yields
\[
\left|\int_0^T\varphi(t)\,dt-\sum_{i=1}^{\delta}\xi_i\varphi(t_i)\right|
\leq C_r\delta^{-(r+1)}
\|\varphi^{(r+1)}\|_{L^\infty(0,T)}.
\]
This bound holds for $\varphi\in W^{r+1,\infty}(0,T)$.
Here $C_r$ depends on $T$, the reference rule, and the mesh regularity. Applying this estimate uniformly requires a common bound on the numerical integrand's derivatives. The control norm alone does not provide that bound.

Let $\mathrm{Quad}(\delta)$ bound
$|\hat{J}^h(u)-\hat{J}^{h,\delta}(u)|$
uniformly for $u\in\mathcal{U}_N$ and the solver meshes under consideration.
Then
\begin{equation}\label{e25}
	\begin{aligned}
		\mathcal{E}_{\mathrm{quad}}
		&=|\hat{J}^h(u_{N,h,\delta}^{k})-\hat{J}^{h,\delta}(u_{N,h,\delta}^{k})|\\
		&\quad+|\hat{J}^{h,\delta}(u_{N,h}^{*})-\hat{J}^h(u_{N,h}^{*})|\\
		&\leq2\,\mathrm{Quad}(\delta),
	\end{aligned}
\end{equation}
where the required consistency condition is
\begin{equation}\label{e26}
	\mathrm{Quad}(\delta)\rightarrow0\qquad(\delta\rightarrow\infty).
\end{equation}

A concrete bound compatible with Lipschitz controls is obtained from the composite midpoint rule. Suppose, in addition, that
\[
|g(x,u,t)-g(x,u,s)|\leq L_{gt}|t-s|,
\qquad \operatorname{Lip}_t(S_h(u))\leq M_f
\]
uniformly in $u$ and $h$, as proved above for the Euler reconstruction. Then
$\varphi_u^h$ is $D_g$-Lipschitz, with
$D_g=L_{gx}M_f+L_{gu}L_0+L_{gt}$.
On $\delta$ equal cells of length $T/\delta$, use their midpoints and weights $\xi_i=T/\delta$. Integrating the bound
$|\varphi_u^h(t)-\varphi_u^h(t_i)|\leq D_g|t-t_i|$ over each cell gives
\[
\sup_{u\in\mathcal{U}_N}
|\hat{J}^h(u)-\hat{J}^{h,\delta}(u)|
\leq\frac{D_gT^2}{4\delta}.
\]
Thus $\mathrm{Quad}(\delta)=D_gT^2/(4\delta)$ verifies \eqref{e26} without higher derivatives of the learned term.

If instead $\varphi_u^h\in W^{2,\infty}$ on each quadrature cell, with a common second-derivative bound $M_2$, the midpoint error on a cell of length $\Delta_i$ is at most $M_2\Delta_i^3/24$. Thus, for midpoint nodes $t_i$ and weights $\Delta_i$,
\[
 \left|\int_0^T\varphi_u^h(t)\,dt
       -\sum_{i=1}^{\delta}\Delta_i\varphi_u^h(t_i)\right|
 \le \frac{M_2}{24}\sum_{i=1}^{\delta}\Delta_i^3
 \le \frac{M_2T}{24}\bigl(\max_i\Delta_i\bigr)^2.
\]
The bound $M_2T^3/(24\delta^2)$ is the equal-cell specialization. Aligning a mesh with spline or solver breakpoints may produce unequal cells, in which case the sum-of-cubes bound applies. Every breakpoint must be a cell boundary, or the crossing cells require a separate estimate. The uniform Lipschitz bound above remains available without this second-derivative assumption.

\subsection{Generalization error}
For data-driven loss functions, let \(Z_i=(t_i,Y_i)\), \(1\leq i\leq\delta\), be independent and identically distributed observations with law \(\mathsf P\), where \(t_i\) is uniform on \([0,T]\). Suppose that the measurable loss satisfies
\[
0\leq\ell(x,a,y)\leq B,\qquad
|\ell(x,a,y)-\ell(\widetilde x,\widetilde a,y)|
\leq L_{gx}|x-\widetilde x|_1+L_{gu}|a-\widetilde a|_\infty
\]
on the exact and numerical state ranges, for admissible controls and all possible observations \(y\). Set
\[
g(x,a,t)=\mathbb E[\ell(x,a,Y)\mid t],\qquad
\hat J(u)=T\mathbb E[\ell(S(u)(t),u(t),Y)].
\]
This expression agrees with the integral objective in \eqref{e1}. Replacing \(S\) by \(S_h\) gives the finite numerical population objective \(\hat J^h\) on \(\mathcal U_N\). For the empirical loss, recall that
\[
\hat J^{h,\delta}(u)=
\begin{cases}
\displaystyle\frac{T}{\delta}\sum_{i=1}^{\delta}
 \ell(S_h(u)(t_i),u(t_i),Y_i),&u\in\mathcal U_N,\\[1ex]
+\infty,&u\in\mathcal U\setminus\mathcal U_N.
\end{cases}
\]
The factor \(T\) preserves the integral normalization; for a mean loss, all objectives and their error bounds are divided by \(T\). In a trajectory observation model, one may take
\(Y_i=\mathcal O S(u^\dagger)(t_i)+\xi_i\), where \(\mathcal O\) is the observation map. The true trajectory is fixed. Taking the pairs \((t_i,\xi_i)\) independently from one common joint law makes \(Z_i\) iid, while allowing the noise distribution to depend on the observation time. This is a random-time observation model, not an independence assumption about successive states of the trajectory. For squared loss, an observation map bounded on the state range and bounded noise ensure the required loss bound on a compact prediction range.

To bound the generalization error \(\mathcal E_{\text{gen}}\), we use Rademacher complexity \cite{bartlett2002rademacher}. A textbook treatment is given in \cite{shalev2014understanding}. The relevant hypothesis class is
\[
\mathcal H_{N,h}
=\{H_{h,u}:(t,y)\mapsto\ell(S_h(u)(t),u(t),y):
                         u\in\mathcal U_N\}.
\]
This definition includes the dependence of the numerical trajectory on the whole control function. The compact parameter representation and Lipschitz assumptions below imply separability in the uniform norm. The suprema can therefore be taken over countable dense subsets and are measurable.

\begin{definition}(Rademacher complexity).
Given a class of measurable functions \(\mathcal H\) and sample points
\(\mathcal S=(Z_1,\ldots,Z_\delta)\), the empirical Rademacher complexity of \(\mathcal H\) is defined by
\[
\hat{\mathcal R}_\delta(\mathcal H)
=\frac1\delta\mathbb E_\sigma
 \left[\sup_{\varphi\in\mathcal H}
       \sum_{i=1}^{\delta}\sigma_i\varphi(Z_i)\right].
\]
\end{definition}
Here \(\sigma_i\), \(1\leq i\leq\delta\), are independent Rademacher random variables, uniformly distributed on \(\{-1,1\}\) and independent of the sample. The Rademacher complexity is the average over the samples:
\[
\mathcal R_\delta(\mathcal H)
=\mathbb E_{Z_1,\ldots,Z_\delta\sim\mathsf P^\delta}
       \hat{\mathcal R}_\delta(\mathcal H).
\]
A standard symmetrization result controls the generalization gap \cite{wainwright2019high}.
\begin{lemma}\label{lemma2}
Let \(\mathcal H\) be a nonempty separable class of measurable functions taking values in \([0,B]\), and put
\(\mathcal H^\pm=\mathcal H\cup\{-\varphi:\varphi\in\mathcal H\}\).
Then
\[
\mathbb E\sup_{\varphi\in\mathcal H}
 \left|\frac1\delta\sum_{i=1}^{\delta}\varphi(Z_i)
                      -\mathbb E\varphi(Z)\right|
\leq2\mathcal R_\delta(\mathcal H^\pm).
\]
For every \(0<\alpha<1\), with probability at least \(1-\alpha\),
\[
\sup_{\varphi\in\mathcal H}
 \left|\frac1\delta\sum_{i=1}^{\delta}\varphi(Z_i)
                      -\mathbb E\varphi(Z)\right|
\leq2\mathcal R_\delta(\mathcal H^\pm)
       +B\sqrt{\frac{\log(1/\alpha)}{2\delta}}.
\]
\end{lemma}
\begin{proof}
Introduce an independent sample with the same law and apply symmetrization. The absolute supremum equals the supremum over \(\mathcal H^\pm\), which gives the expectation bound. Changing one observation changes the absolute empirical supremum by at most \(B/\delta\). The bounded-differences inequality therefore gives the stated probability bound.
\end{proof}

Here the Rademacher complexity depends on the entire loss class, and its decay must be justified for that class. The scalar contraction inequality in \cite{shalev2014understanding} concerns pointwise composition with a Lipschitz function. It does not apply directly to the map that takes a control function to its ODE trajectory. We give a finite-parameter estimate for \(\mathcal H_{N,h}\) in the following proposition.
\begin{proposition}[A covering bound for the trajectory loss class]\label{prop:loss_covering}
Suppose \(\mathcal U_N=\{u_\theta:\theta\in\Theta_N\}\), where \(\Theta_N\) is a nonempty compact subset of a \(P_N\)-dimensional coordinate box with side lengths at most \(D_N\), and
\[
\|u_\theta-u_\vartheta\|_\infty
\leq A_N|\theta-\vartheta|_\infty.
\]
Set \(L_H=(L_{gx}L_{S,h}+L_{gu})A_N\). For each fixed \(r>0\) and \(0<\alpha<1\), define
\begin{equation}\label{eq:loss_covering}
\begin{aligned}
\zeta_{\delta,N,h}(r,\alpha)
&:=2L_Hr\\
&\quad+B\left(
 \frac{P_N\log(1+D_N/r)+\log(2/\alpha)}{2\delta}
                  \right)^{1/2}.
\end{aligned}
\end{equation}
Then, with probability at least \(1-\alpha\),
\[
\sup_{u\in\mathcal U_N}
 |\hat J^{h,\delta}(u)-\hat J^h(u)|
\leq T\zeta_{\delta,N,h}(r,\alpha).
\]
\end{proposition}
\begin{proof}
The loss and state-map bounds give
\[
\|H_{h,u_\theta}-H_{h,u_\vartheta}\|_\infty
\leq L_H|\theta-\vartheta|_\infty.
\]
Partition the parameter box into coordinate cells of diameter at most \(r\) in the maximum norm. Choose one point of \(\Theta_N\) in each cell that intersects \(\Theta_N\). This gives an admissible \(r\)-net with at most \((1+D_N/r)^{P_N}\) points. Hoeffding's inequality \cite{hoeffding1963probability} and a union bound imply that, with probability at least \(1-\alpha\), the normalized empirical deviation at every net point is at most
\[
B\left(
 \frac{P_N\log(1+D_N/r)+\log(2/\alpha)}{2\delta}
 \right)^{1/2}.
\]
Both the empirical and population losses change by at most \(L_Hr\) between an arbitrary parameter and its net point. Adding these two errors and restoring the factor \(T\) proves the result. If \(B=0\), every loss is zero and the assertion is immediate.
\end{proof}

For the nodal ReLU class, one can take \(P_N=m(N+1)\), \(A_N=1\), and \(D_N=\max_i(b_i-a_i)\). If \(L_{S,h}\) is bounded uniformly over the meshes used, choosing \(r=\delta^{-1}\) gives convergence of this bound whenever
\((N+1)\log\delta/\delta\to0\) at fixed confidence. This is a sufficient condition for the displayed parameter estimate. Section~\ref{sec3} also gives qualitative uniform convergence on the fixed compact control class.

The event in Proposition~\ref{prop:loss_covering} is uniform over controls in the fixed class, not automatically over resolutions selected from the same data. For a finite or countable family of choices specified in advance $(N_j,h_j,\delta_j,r_j)$, apply the proposition with levels $\alpha_j>0$ such that $\sum_j\alpha_j\le\alpha$. A union bound then makes all corresponding estimates valid simultaneously with probability at least $1-\alpha$, allowing selection within this family after observing the data.

We can now bound the statistical term in Theorem~\ref{thm:error_bound}. On the event in Proposition \ref{prop:loss_covering}, both controls in the definition of \(\mathcal E_{\text{gen}}\) belong to \(\mathcal U_N\), and hence
\begin{equation}\label{e23}
\begin{aligned}
\mathcal E_{\text{gen}}
&=|\hat J^h(u_{N,h,\delta}^k)
                  -\hat J^{h,\delta}(u_{N,h,\delta}^k)|\\
&\quad+|\hat J^{h,\delta}(u_{N,h}^*)-\hat J^h(u_{N,h}^*)|\\
&\leq2\sup_{u\in\mathcal U_N}
                     |\hat J^{h,\delta}(u)-\hat J^h(u)|\\
&\leq2T\zeta_{\delta,N,h}(r,\alpha).
\end{aligned}
\end{equation}
The same event covers every feasible iterate, including one chosen using the observations. Alternatively, Lemma~\ref{lemma2} gives
\[
\mathcal E_{\text{gen}}
\leq4T\mathcal R_\delta(\mathcal H_{N,h}^\pm)
       +2TB\sqrt{\frac{\log(1/\alpha)}{2\delta}}
\]
with probability at least \(1-\alpha\). Either estimate can be used in Theorem \ref{thm:error_bound}.

\subsection{Optimization error}
Let \(u_{N,h,\delta}^k\in\mathcal U_N\) be the feasible control returned after \(k\) optimization steps. We measure its optimization error by the objective gap
\begin{equation}\label{e19}
\begin{aligned}
\rho_{N,h,\delta}(k)
&:=\hat J^{h,\delta}(u_{N,h,\delta}^k)
             -\min_{u\in\mathcal U_N}\hat J^{h,\delta}(u)
\geq0.
\end{aligned}
\end{equation}
If a specified optimization method admits an upper bound
\(\rho_{N,h,\delta}(k)\leq C_{\mathrm{opt}}\mathrm{Opt}(k)\),
we may insert that bound into the error decomposition. Under this additional condition,
\begin{equation}\label{e20}
\begin{aligned}
\mathcal E_{\text{opt}}
&=|\hat J^{h,\delta}(u_{N,h,\delta}^k)
                   -\hat J^{h,\delta}(u_{N,h,\delta}^*)|\\
&=\rho_{N,h,\delta}(k)
\leq C_{\mathrm{opt}}\mathrm{Opt}(k).
\end{aligned}
\end{equation}
The identity with \(\rho_{N,h,\delta}(k)\) holds for every feasible iterate. The upper bound and the limit \(\mathrm{Opt}(k)\to0\) require an optimization result for the problem under consideration.

The convergence rate depends on the convexity and smoothness of the parameter objective and on the optimization algorithm. The following classical rates concern objective gaps under their stated assumptions; constants can depend on the fixed discretized problem.
For deterministic first-order methods applied to smooth convex objectives on a convex parameter domain, using the corresponding projected methods when constraints are present:
\begin{itemize}
\item Gradient Descent (GD): A suitable step size gives an \(O(k^{-1})\) objective gap \cite{nesterov2013introductory}. If the objective is strongly convex, the convergence is linear, with an \(O(q^k)\) gap for some \(0<q<1\).
\item Accelerated Gradient Methods (e.g., Nesterov's acceleration): These methods attain an \(O(k^{-2})\) objective gap for smooth convex objectives \cite{nesterov2013introductory}.
\end{itemize}
For stochastic first-order methods, the gradient is estimated from mini-batches, and the analysis must also control the variance of that estimate.
\begin{itemize}
\item Stochastic Gradient Descent (SGD): For an unconstrained, lower-bounded smooth nonconvex objective with unbiased gradients of bounded variance, suitable step sizes give an \(O(k^{-1/2})\) bound on the smallest expected squared gradient norm among the first \(k\) iterates \cite{bottou2018optimization}. This is a stationarity estimate and does not bound the global objective gap in \eqref{e19}. Adaptive methods such as Adam \cite{kingma2014adam} require their own assumptions; vanilla Adam can fail to converge even on convex examples \cite{reddi2019convergence}.
\item SGD for Convex Problems: Under strong convexity and a Lipschitz gradient, unbiased gradient estimates with bounded variance and suitable diminishing step sizes give \(O(k^{-1})\) expected objective gaps \cite{bottou2018optimization}. Here the expectation is over the optimization randomness for a fixed training sample.
\end{itemize}
Global minimization of a nonconvex objective requires stronger algorithm-specific assumptions. Some Bayesian optimization schemes based on Gaussian processes need not explore the design space densely; expected improvement can be inconsistent even for smooth objectives \cite{yarotsky2013examples}. Thus a vanishing global objective gap requires assumptions on the particular optimization procedure.

In the context of deep learning for UDEs, the objective \(\hat J^{h,\delta}\) is generally nonconvex, and results about stationary points do not identify a global minimizer. Theorem \ref{thm:error_bound} therefore keeps \(\rho_{N,h,\delta}(k)\) as an explicit error term. When the discretization changes with \(j\), consistency of the returned controls requires a choice of iterates for which
\(\rho_{N_j,h_j,\delta_j}(k_j)\to0\); the condition \(k_j\to\infty\) alone does not supply this conclusion.

\section{Proof of \texorpdfstring{$\Gamma$}{Gamma}-Convergence Results}\label{sec3}
We analyze the approximate objective functionals \(\{\hat J_\epsilon\}\) through \(\Gamma\)-convergence \cite{braides2002gamma,brezis1993introduction}. This complements the error estimates in Section \ref{sec2} by treating the changing neural network classes as part of the variational problem.

\begin{definition}(\(\Gamma\)-convergence, \cite{braides2002gamma,brezis1993introduction})
Let \(X\) be a metric space and let \(F_\epsilon,F:X\to\mathbb R\cup\{+\infty\}\). Then \(F_\epsilon\) is said to \(\Gamma\)-converge to \(F\) if the following properties hold for every \(x\in X\).

\noindent (i) Liminf inequality: for every sequence \(x_\epsilon\to x\),
\[
F(x)\leq\liminf_{\epsilon\to0}F_\epsilon(x_\epsilon).
\]

\noindent (ii) Recovery sequence: there exists a sequence \(x_\epsilon\to x\) such that
\[
F(x)=\lim_{\epsilon\to0}F_\epsilon(x_\epsilon).
\]
\end{definition}
There are equivalent definitions of \(\Gamma\)-convergence; we refer the reader to \cite{braides2002gamma} for details.

\begin{definition}(Equi-coercivity, \cite{braides2002gamma,brezis1993introduction})
A sequence \(F_\epsilon:X\to\mathbb R\cup\{+\infty\}\) is equi-coercive if, for every \(c\in\mathbb R\), there exists a compact set \(K_c\subset X\) such that
\(\{x\in X:F_\epsilon(x)\leq c\}\subset K_c\) for every \(\epsilon\).
\end{definition}

\begin{theorem}(Fundamental theorem of \(\Gamma\)-convergence, \cite{braides2002gamma,brezis1993introduction})\label{t3}
Let \((X,d)\) be a metric space, and let \(F_\epsilon,F:X\to\mathbb R\cup\{+\infty\}\) be proper lower semicontinuous functionals. Suppose \(F_\epsilon\) is equi-coercive and \(\Gamma\)-converges to \(F\). Then the minima exist, are finite, and
\[
\min_X F=\lim_{\epsilon\to0}\min_X F_\epsilon.
\]
Moreover, if
\[
F_\epsilon(x_\epsilon)\leq\min_XF_\epsilon+\rho_\epsilon,
\qquad \rho_\epsilon\geq0,\qquad \rho_\epsilon\to0,
\]
then \(x_\epsilon\) is precompact and every limit point is a minimizer of \(F\).
Here proper means that the functional is finite at at least one point.
\end{theorem}
We apply this result to the functionals in Theorem~\ref{thm:gamma_conv}.

\subsection{Assumptions and Setup}
\label{subsec:assumptions}
Let \((\U,d)\) be the metric space of admissible control functions specified in Section \ref{sec2}, with \(d(u,v)=\|u-v\|_\infty\). The exact objective functional is
\[
\hat J(u)=J(S(u),u),\qquad u\in\U,
\]
where \(S:\U\to\X\) is the control-to-state map for the continuous dynamical system. The approximate objective depends on the ODE step size \(h\), the number of quadrature nodes or observations \(\delta\), and the network complexity \(N\). Choose a sequence of triples satisfying
\[
h_j\to0,\qquad \delta_j\to\infty,\qquad N_j\to\infty,
\]
and define
\[
\epsilon_j=\max\left\{h_j,\frac1{\delta_j},\frac1{N_j}\right\},
\qquad
\hat J_{\epsilon_j}=\hat J^{h_j,\delta_j}.
\]
The notation \(\epsilon_j\) indexes the chosen sequence of triples. The condition \(\epsilon_j\to0\) is equivalent to the three limits above. The dependence on \(N_j\) enters through the domain:
\[
\hat J_{\epsilon_j}(u)=+\infty\quad\hbox{for }u\in\U\setminus\U_{N_j}.
\]
For empirical objectives, the exact functional is the population integral defined in the generalization-error subsection.

To establish \(\Gamma\)-convergence, we use the following assumptions.
\begin{assumption}[Compactness of \(\U\)]
\label{ass:compact}
The space \((\U,d)\) is a compact metric space. In the present setting,
\[
\U=\{u\in C([0,T];A):\operatorname{Lip}_\infty(u)\leq L_0\},
\]
where \(A=\prod_{i=1}^m[a_i,b_i]\), and \(L_0\) is fixed throughout refinement. This class is uniformly bounded and equicontinuous. Its range and Lipschitz constraints are closed under uniform convergence, so compactness follows from the Arzel\`a--Ascoli theorem.
\end{assumption}

\begin{assumption}[Continuity of \(\hat J\)]
\label{ass:cont}
The exact objective functional \(\hat J:\U\to\mathbb R\) is continuous with respect to \(d\). In particular, Lemma \ref{lemma1} gives its Lipschitz constant \(T(L_{gx}L_S+L_{gu})\).
\end{assumption}

\begin{assumption}[Consistency of numerical solver]
\label{ass:solver}
For every mesh under consideration, the numerical state map \(S_h:\U\to\X\) is continuous, and its uniform absolute error satisfies
\[
\mathrm{Num}(h)
:=\sup_{u\in\U}\|S_h(u)-S(u)\|_{\mathcal X}\longrightarrow0
\qquad(h\to0).
\]
The exact and numerical trajectories remain in the common domain on which the cost bounds hold.
\end{assumption}

\begin{assumption}[Density of the neural network classes]
\label{ass:density}
Each neural network class \(\U_N\) is a nonempty compact subset of \(\U\), and
\[
\operatorname{dist}(u,\U_N)\longrightarrow0
\qquad(N\to\infty)
\quad\hbox{for every }u\in\U.
\]
The finite restriction of \(\hat J^{h,\delta}\) to \(\U_N\) is continuous. The classes need not be nested.
\end{assumption}

\begin{assumption}[Uniform error bound]
\label{ass:uniform}
For the chosen refinement sequence, there are numbers \(\eta_j\geq0\) such that \(\eta_j\to0\) and
\[
\sup_{u\in\U_{N_j}}
 |\hat J^{h_j,\delta_j}(u)-\hat J(u)|\leq\eta_j.
\]
The comparison is restricted to \(\U_{N_j}\), where both objectives are finite. For random objectives, this assumption and the resulting variational conclusions are interpreted on the sample event specified below.
\end{assumption}

These assumptions are compatible with the approximation and numerical estimates in Section \ref{sec2}. The nodal ReLU construction gives a compact parameter image contained in \(\U\), with
\(\operatorname{dist}(u,\U_N)\leq L_0T/(2N)\).
Continuity of \(S_h\), point evaluation in the uniform norm, and the loss bounds give continuity of the finite numerical objectives. Compactness of \(\U_N\) therefore ensures existence of their minimizers. Since \(\U_N\) is closed in \(\U\), extension by \(+\infty\) yields a proper lower semicontinuous functional. Assumption \ref{ass:solver} is used to derive the finite-domain error in Assumption \ref{ass:uniform}, as follows.

\begin{lemma}[Derivation of Assumption \ref{ass:uniform} from Section \ref{sec2}]
For deterministic quadrature, the uniform quadrature estimate gives
\[
\sup_{u\in\U_N}
 |\hat J^{h,\delta}(u)-\hat J(u)|
\leq TL_{gx}\mathrm{Num}(h)+\mathrm{Quad}(\delta).
\]
For the empirical objective, under the assumptions of Proposition \ref{prop:loss_covering}, with probability at least \(1-\alpha\),
\[
\sup_{u\in\U_N}
 |\hat J^{h,\delta}(u)-\hat J(u)|
\leq TL_{gx}\mathrm{Num}(h)+T\zeta_{\delta,N,h}(r,\alpha).
\]
The first estimate applies to deterministic quadrature; the second applies to random observations. If the corresponding bound tends to zero along the chosen sequence, it establishes Assumption~\ref{ass:uniform} on the event where that bound holds.
\end{lemma}
\begin{proof}
For \(u\in\U_N\), the error decomposition in Section \ref{sec2} gives
\[
|\hat J^{h,\delta}(u)-\hat J(u)|
\leq|\hat J^h(u)-\hat J(u)|
     +|\hat J^{h,\delta}(u)-\hat J^h(u)|.
\]
The first term is at most \(TL_{gx}\mathrm{Num}(h)\) by integration of the state error and the Lipschitz bound for \(g\). The second term is bounded by the uniform quadrature error in \eqref{e25}--\eqref{e26}, or by Proposition \ref{prop:loss_covering} for the empirical loss. Taking the supremum over \(\U_N\) proves the claims.
\end{proof}

For statistical confidence bounds along a deterministic refinement sequence, choose deterministic \(r_j>0\) and \(0<\alpha_j<1\) such that
\[
\sum_j\alpha_j<\infty,\qquad
\zeta_{\delta_j,N_j,h_j}(r_j,\alpha_j)\longrightarrow0.
\]
The Borel--Cantelli lemma implies that the empirical bound holds for all sufficiently large \(j\), almost surely. No independence between samples used at different refinement levels is needed for this implication.

There is also a qualitative statistical conclusion on the fixed compact class \(\U\). Take observations from prefixes of one infinite iid sequence \(Z_1,Z_2,\ldots\), and consider the exact trajectory loss class
\[
\mathcal H
=\{H_u:(t,y)\mapsto\ell(S(u)(t),u(t),y):u\in\U\}.
\]
Its elements are bounded by \(B\), and
\[
\|H_u-H_v\|_\infty
\leq (L_{gx}L_S+L_{gu})\|u-v\|_\infty.
\]
Thus \(\mathcal H\) is compact in the uniform norm. For each positive integer \(q\), choose a finite \(1/q\)-net in \(\mathcal H\). The strong law holds at every net point, along all sample prefixes. Intersecting these countably many probability-one events yields a single event on which
\[
\begin{aligned}
\Delta_\delta
&:=\sup_{u\in\U}
 \left|\frac1\delta\sum_{i=1}^{\delta}H_u(Z_i)-\mathbb EH_u(Z)\right|,\\
\limsup_{\delta\to\infty}\Delta_\delta&\leq\frac2q
\qquad\hbox{for every }q.
\end{aligned}
\]
Letting \(q\to\infty\) gives \(\Delta_\delta\to0\) almost surely. The effect of the numerical state error on the empirical loss is uniform over the same controls, so
\[
\sup_{u\in\U_{N_j}}
 |\hat J^{h_j,\delta_j}(u)-\hat J(u)|
\leq TL_{gx}\mathrm{Num}(h_j)+T\Delta_{\delta_j}
\longrightarrow0
\]
on this event. Hence Assumption \ref{ass:uniform} holds almost surely for every refinement sequence with \(h_j\to0\), \(\delta_j\to\infty\), and \(N_j\to\infty\), when data are taken from these prefixes. This qualitative result uses the fixed bound \(L_0\); it imposes no additional coupling of \(N_j\) and \(\delta_j\). Formula \eqref{eq:loss_covering} gives a separate sufficient condition when an explicit confidence bound is needed.

\subsection{Proof of \texorpdfstring{$\Gamma$}{Gamma}-Convergence Result}
\label{subsec:gamma-result}
\begin{proof}
(1) We divide the proof into three steps: equi-coercivity, the liminf inequality, and the construction of a recovery sequence. For empirical objectives, fix a sample realization in an event where Assumption \ref{ass:uniform} holds, such as the probability-one event established above.

\noindent Step 1: Equi-coercivity.
By Assumption \ref{ass:compact}, the whole space \(\U\) is compact. Hence, for any \(c\in\mathbb R\), the sublevel set
\(\{u\in\U:\hat J_{\epsilon_j}(u)\leq c\}\)
is contained in the compact set \(\U\). Therefore, the sequence \(\{\hat J_{\epsilon_j}\}\) is equi-coercive.

\noindent Step 2: Liminf inequality.
Let \(u_j\to u\) in \(\U\). We need to show
\[
\hat J(u)\leq\liminf_{j\to\infty}\hat J_{\epsilon_j}(u_j).
\]
If the lower limit is \(+\infty\), the inequality holds trivially. Otherwise, take a subsequence that realizes it and has finite objective values. These values are bounded below by \(\min_\U\hat J-\eta_j\), so the lower limit cannot be \(-\infty\). Along this subsequence, \(u_j\in\U_{N_j}\). Assumption \ref{ass:uniform} then gives
\[
\hat J_{\epsilon_j}(u_j)\geq\hat J(u_j)-\eta_j.
\]
By Assumption \ref{ass:cont}, \(\hat J(u_j)\to\hat J(u)\), while \(\eta_j\to0\). Taking the lower limit yields
\[
\liminf_{j\to\infty}\hat J_{\epsilon_j}(u_j)\geq\hat J(u).
\]

\noindent Step 3: Recovery sequence.
For any \(u\in\U\), we construct \(v_j\to u\) such that
\[
\lim_{j\to\infty}\hat J_{\epsilon_j}(v_j)=\hat J(u).
\]
By compactness of \(\U_{N_j}\) and Assumption \ref{ass:density}, there is \(v_j\in\U_{N_j}\) with
\[
d(u,v_j)=\operatorname{dist}(u,\U_{N_j})\longrightarrow0.
\]
Using Assumption \ref{ass:uniform} again,
\[
\begin{aligned}
|\hat J_{\epsilon_j}(v_j)-\hat J(u)|
&\leq|\hat J_{\epsilon_j}(v_j)-\hat J(v_j)|
                      +|\hat J(v_j)-\hat J(u)|\\
&\leq\eta_j+|\hat J(v_j)-\hat J(u)|.
\end{aligned}
\]
The first term vanishes by Assumption \ref{ass:uniform}, and the second vanishes by continuity of \(\hat J\). Hence \(v_j\) is a recovery sequence.

Therefore, \(\{\hat J_{\epsilon_j}\}\) is equi-coercive and \(\Gamma\)-converges to \(\hat J\). Both the exact and approximate functionals are proper and lower semicontinuous, with minima attained on their compact feasible sets. Applying Theorem \ref{t3} yields
\[
\min_{u\in\U}\hat J(u)
=\lim_{j\to\infty}\min_{u\in\U}\hat J_{\epsilon_j}(u),
\]
which completes the proof of (1).

\noindent (2) Let \(\{u_{\epsilon_j}^*\}\) be any sequence such that
\[
u_{\epsilon_j}^*\in\arg\min_{u\in\U}\hat J_{\epsilon_j}(u).
\]
Since \(\U\) is compact, the sequence is precompact: every subsequence has a further subsequence converging in \(\U\). Let \(u\) be such a limit point. Write \(m=\min_\U\hat J\), and use the liminf inequality together with the convergence of minimum values from (1):
\[
\hat J(u)\leq\liminf_{j\to\infty}
                   \hat J_{\epsilon_j}(u_{\epsilon_j}^*)=m,
\]
where the indices are restricted to the convergent subsequence. Since \(m\) is the minimum, \(\hat J(u)=m\), so \(u\in\arg\min_\U\hat J\). This proves the asserted convergence of minimizers.

The same argument applies to feasible approximate minimizers \(u_j\in\U_{N_j}\) satisfying
\[
\hat J_{\epsilon_j}(u_j)
\leq\min_\U\hat J_{\epsilon_j}+\rho_j,\qquad
\rho_j\geq0,\qquad \rho_j\to0.
\]
Indeed, their approximate objective values converge to \(m\). Assumption~\ref{ass:uniform} also gives \(\hat J(u_j)\to m\), and every limit point belongs to the nonempty compact set $\mathcal M=\arg\min_{u\in\mathcal U}\hat J(u)$. In particular,
\[
 \operatorname{dist}(u_j,\mathcal M)
 :=\inf_{v\in\mathcal M}\|u_j-v\|_\infty\longrightarrow0.
\]
Otherwise, some subsequence would stay a positive distance from $\mathcal M$. Compactness would give a further convergent subsequence whose limit belongs to $\mathcal M$, a contradiction. If the minimizer is unique, this is convergence of the whole sequence to that control.
\end{proof}

\begin{remark}
\label{rmk:interpretation}
Theorem \ref{thm:gamma_conv} guarantees that the minima of the approximate problems converge to the minimum of the continuous problem. Moreover, any sequence of minimizers of \(\hat J_{\epsilon_j}\) admits a subsequence converging to a minimizer of \(\hat J\). For controls returned by training, this conclusion requires the global objective gaps in \eqref{e19} to tend to zero. Without uniqueness, convergence is to the set of minimizers and does not identify a particular underlying control.
\end{remark}
\begin{remark}
If \(\U\) is not compact, one must establish equi-coercivity and continuity or lower semicontinuity in a specified topology. The compactness of \(\U\) simplifies the present argument. In this setting, convergence of minimum values can also be proved by elementary comparison estimates. The \(\Gamma\)-convergence formulation additionally records the limiting behavior of the changing feasible classes and their extensions by \(+\infty\); the uniform finite-domain bound in Assumption \ref{ass:uniform} is not uniform convergence of the extended functionals on all of \(\U\).
\end{remark}

\section{Example: Learning a Time-Varying Transmission Rate in the SIR Model}
\label{sec:sir_example}

We apply the preceding estimates to the Susceptible-Infected-Recovered
(SIR) model with an unknown time-varying transmission rate. The known
compartmental equations determine how errors in the learned coefficient
affect the observed infected population. The example gives explicit expressions for the four error
terms and examines them in a synthetic inverse problem.

\subsection{Model Setup}
\label{subsec:sir_setup}

Consider the SIR model with time-varying transmission rate \(\beta(t)\) and
constant recovery rate \(\gamma>0\):
\begin{align}
 &\frac{\mathrm{d}S}{\mathrm{d}t}=-\beta(t)SI,
     \label{eq:sir_beta_s}\\
 &\frac{\mathrm{d}I}{\mathrm{d}t}=\beta(t)SI-\gamma I,
     \label{eq:sir_beta_i}\\
 &\frac{\mathrm{d}R}{\mathrm{d}t}=\gamma I.
     \label{eq:sir_beta_r}
\end{align}
The initial state \(x_0=(S_0,I_0,R_0)\) is known, nonnegative, and has unit
total population. The simplex
\(K=\{(S,I,R)\geq0:S+I+R=1\}\) is invariant: the vector field points inward
on each coordinate boundary, and its components sum to zero. We take the
admissible class \(\mathcal B_{L_0}\) to consist of continuous controls
satisfying \(\beta_{\min}\leq\beta(t)\leq\beta_{\max}\) and
\(\operatorname{Lip}(\beta)\leq L_0\), where
\(0\leq\beta_{\min}<\beta_{\max}\). This is the scalar instance of the compact
control class used above.

The transmission rate is learned from observations of \(I\). For the
statistical estimates, let \(\beta^*\in\mathcal B_{L_0}\) generate the data
and assume that the pairs \((t_i,\xi_i)\) are independent and identically
distributed, with
\[
 t_i\sim\operatorname{Unif}[0,T],\qquad
 Y_i=I_{\beta^*}(t_i)+\xi_i,\qquad
 \mathbb E[\xi_i\mid t_i]=0,\quad |\xi_i|\leq\sigma.
\]
Thus independence is across observation pairs; it is not an assumption
about successive states of a deterministic trajectory. We
parameterize \(\beta\) by a feasible network \(\beta_\theta(t)\), giving
\[
 \frac{\mathrm{d}x}{\mathrm{d}t}=f(x,\beta_\theta(t),t),
 \qquad x=(S,I,R)^\top,
\]
where \(f\) is specified by
\eqref{eq:sir_beta_s}--\eqref{eq:sir_beta_r} and has no explicit time
dependence apart from the control.

Let \(I_\beta^h\) denote the infected component of a continuous numerical
trajectory. The empirical mean squared loss is
\begin{equation}\label{eq:sir_empirical_loss}
 \mathcal L(\theta)=\mathcal L_{M,h}(\beta_\theta)
 =\frac1M\sum_{i=1}^M\bigl(I_{\beta_\theta}^h(t_i)-Y_i\bigr)^2.
\end{equation}
To distinguish training error from prediction error, define the population
losses as
\[
 \mathcal L_{\rm pop}(\beta)=\mathbb E(I_\beta(t)-Y)^2,
 \qquad
 \mathcal L_h(\beta)=\mathbb E(I_\beta^h(t)-Y)^2.
\]
These losses are normalized averages, whereas the time-integral objective
\(\hat J\) in the preceding sections is unnormalized. With
\(\nu_\xi=\mathbb E\xi^2\), conditional centering gives
\begin{equation}\label{eq:sir_population_identity}
 \mathcal L_{\rm pop}(\beta)-\nu_\xi
 =\frac1T\int_0^T\bigl(I_\beta(t)-I_{\beta^*}(t)\bigr)^2\,\mathrm dt.
\end{equation}
Hence \(\mathcal L^*:=\min_{\mathcal B_{L_0}}\mathcal L_{\rm pop}=\nu_\xi\).
The explicit numerical estimates below use Euler integration. Adaptive
solvers such as Tsit5 may also be used to evaluate the dynamics, but their
error tolerances require a separate analysis before entering an objective
error bound. We similarly distinguish the choice of optimizer from a
guarantee on its final objective gap.

\subsection{Error Decomposition with Explicit Constants}
\label{subsec:sir_explicit_errors}

Let \(\widehat\beta=\beta_{\theta_K}\) be a feasible control returned after
\(K\) optimization iterations. The four-term comparison now concerns the
population objective:
\[
 \mathcal L_{\rm pop}(\widehat\beta)-\mathcal L^*
 \leq\mathcal E_{\rm approx}+\mathcal E_{\rm num}
       +\mathcal E_{\rm gen}+\mathcal E_{\rm opt}.
\]
We specify each term for an admissible ReLU class and a numerical scheme
that preserve the bounds on the control and state.

\paragraph{1. Approximation Error}

ReLU approximation estimates depend on the regularity of the target and
the network architecture \cite{yarotsky2017error}. Here an explicit
one-dimensional construction suffices. Let \(\mathcal U_N\) be the class of
continuous linear splines on \(N\) equal intervals, written as
one-hidden-layer ReLU networks. Their nodal values \(c_0,\ldots,c_N\) obey
\[
 \beta_{\min}\leq c_j\leq\beta_{\max},\qquad
 |c_{j+1}-c_j|\leq L_0T/N.
\]
The linear interpolant of every \(\beta\in\mathcal B_{L_0}\) satisfies these
constraints. The nodal coefficient set is compact, so a best approximation
exists. Its uniform error satisfies
\begin{equation}\label{eq:sir_approximation_rate}
 d_N:=\inf_{\beta\in\mathcal U_N}\|\beta-\beta^*\|_\infty
 \leq\frac{L_0T}{2N}.
\end{equation}
If \(\beta^*\in W^{2,\infty}(0,T)\), the same interpolant also gives
\(d_N\leq T^2\|\beta^{*\prime\prime}\|_\infty/(8N^2)\). This additional
regularity is what permits a second-order coefficient approximation.

To bound the state sensitivity, we estimate both the state Jacobian and the
derivative with respect to the transmission rate. For the SIR equations,
\[
 D_xf=
 \begin{pmatrix}
 -\beta I&-\beta S&0\\
 \beta I&\beta S-\gamma&0\\
 0&\gamma&0
 \end{pmatrix},
 \qquad
 \partial_\beta f=(-SI,SI,0)^\top.
\]
On the convex simplex \(K\),
\begin{align*}
 \|D_xf\|_1
 &=\max\{2\beta I,\,\beta S+|\beta S-\gamma|+\gamma\}
   \leq2\max\{\beta_{\max},\gamma\},\\
 \|\partial_\beta f\|_1&=2SI\leq\tfrac12.
\end{align*}
We may therefore take
\begin{equation}\label{eq:sir_sensitivity_constants}
 L_{fx}=2\max\{\beta_{\max},\gamma\},\qquad
 L_{fu}=\tfrac12,\qquad L_S=\tfrac12\Phi(L_{fx},T).
\end{equation}
The control-to-state estimate then reads
\(\|x_\beta-x_{\widetilde\beta}\|_{\mathcal X}
\leq L_S\|\beta-\widetilde\beta\|_\infty\), with the uniform-in-time
\(\ell^1\) state norm. Combining it with
\eqref{eq:sir_population_identity} yields
\[
 \mathcal E_{\rm approx}
 :=\min_{\mathcal U_N}\mathcal L_{\rm pop}-\mathcal L^*
 \leq L_S^2d_N^2.
\]
The square follows from the centered-noise identity, not from a general
Lipschitz estimate for an arbitrary objective.

\paragraph{2. Numerical Error}

For a uniform Euler grid with step \(h=T/J\), the updates are
\begin{align*}
 S_{j+1}&=S_j(1-h\beta_j I_j),\\
 I_{j+1}&=I_j(1-h\gamma+h\beta_jS_j),\qquad
 R_{j+1}=R_j+h\gamma I_j.
\end{align*}
If \(h\max\{\beta_{\max},\gamma\}\leq1\), all components remain
nonnegative and their sum remains one. Linear reconstruction between
successive states also lies in \(K\). Moreover,
\[
 \|f\|_1=2\max\{\beta SI,\gamma I\}
 \leq M_f:=\max\{\beta_{\max}/2,2\gamma\}.
\]
The numerical estimate in Section~\ref{sec2} applies with
\(L_{ft}=0\) and \(D_f=L_{fx}M_f+L_0/2\), giving
\begin{equation}\label{eq:sir_euler_bound}
 \sup_{\beta\in\mathcal B_{L_0}}
 \|x_\beta-x_\beta^h\|_{\mathcal X}
 \leq E_h:=D_f\left\{\frac h2\Phi(L_{fx},T)+\frac{h^2}{8}\right\}.
\end{equation}
Both exact and numerical infected fractions lie in \([0,1]\), while
\(Y\in[-\sigma,1+\sigma]\). Thus the squared loss is
\(2(1+\sigma)\)-Lipschitz in the predicted infected fraction, and
\[
 \sup_{\mathcal B_{L_0}}|\mathcal L_{\rm pop}-\mathcal L_h|
 \leq2(1+\sigma)E_h,\qquad
 \mathcal E_{\rm num}\leq4(1+\sigma)E_h.
\]
The two contributions arise from evaluating the learned control and a
comparison control. This is a first-order state bound; it makes no assertion
that an adaptive tolerance equals a global state error.

\paragraph{3. Generalization Error}

The relevant hypothesis class consists of the composed losses
\((t,y)\mapsto(I_\beta^h(t)-y)^2\), \(\beta\in\mathcal U_N\).
General network estimates, including the spectral margin bounds in
\cite{bartlett2017spectrally}, require their specific loss and complexity
assumptions. For the present class we use its compact nodal parameterization.
The Euler difference recursion permits the choice \(L_{S,h}=L_S\) as a
state-map Lipschitz constant. Since linear interpolation is nonexpansive in the
nodal maximum norm, the loss class has parameter Lipschitz constant
\(L_H=2(1+\sigma)L_{S,h}\) and envelope \(B=(1+\sigma)^2\).

Write \(D=\beta_{\max}-\beta_{\min}\). Covering the admissible nodal vectors
by cubes of maximum-norm diameter \(r>0\) requires at most
\((1+D/r)^{N+1}\) representatives, chosen within the feasible set.
For deterministic \(r>0\) and \(0<\alpha<1\), set
\begin{equation}\label{eq:sir_uniform_sampling}
 \begin{aligned}
 \eta_{M,N}(r,\alpha)&:=2L_Hr\\
 &\quad+B\sqrt{\frac{(N+1)\log(1+D/r)+\log(2/\alpha)}{2M}}.
 \end{aligned}
\end{equation}
The finite-cover estimate in Section~\ref{sec2} gives
\[
 \sup_{\beta\in\mathcal U_N}
 |\mathcal L_h(\beta)-\mathcal L_{M,h}(\beta)|
 \leq\eta_{M,N}(r,\alpha)
\]
with probability at least \(1-\alpha\). Hence
\(\mathcal E_{\rm gen}\leq2\eta_{M,N}(r,\alpha)\) on this event. This
uniform statement applies to a control selected using the same observations.
Deterministic quadrature is a separate observation scheme and uses the
quadrature estimate instead of \eqref{eq:sir_uniform_sampling}.

\paragraph{4. Optimization Error}

The optimization term is the global empirical objective gap
\begin{equation}\label{eq:sir_optimization_gap}
 \mathcal E_{\rm opt}=\rho_K
 :=\mathcal L_{M,h}(\widehat\beta)
   -\min_{\beta\in\mathcal U_N}\mathcal L_{M,h}(\beta)\geq0.
\end{equation}
Its decay is an additional condition on the computed controls. It does not
follow from the choice of Adam or from a small gradient. Reddi et al.
\cite{reddi2019convergence} give a convex example in which Adam fails to
converge to the optimum. The nonconvex analysis in
\cite{he2023convergence} imposes restrictions on hyperparameters and, for
function-value rates, a Polyak--\L{}ojasiewicz condition. Such assumptions
must be verified for the objective under consideration. A stationarity
estimate alone does not bound the global gap in
\eqref{eq:sir_optimization_gap}. We therefore retain \(\rho_K\) explicitly
and examine an independent bound in the scalar numerical example below.

\subsection{Total Error Bound}
\label{subsec:sir_total_bound}

To connect the squared approximation estimate to the four-term comparison, set
\[
 \varepsilon_{\rm num}=\sup_{\beta\in\mathcal U_N}
       |\mathcal L_{\rm pop}(\beta)-\mathcal L_h(\beta)|,
 \qquad
 \varepsilon_{\rm samp}=\sup_{\beta\in\mathcal U_N}
       |\mathcal L_h(\beta)-\mathcal L_{M,h}(\beta)|.
\]
Let $\beta_N^{\rm opt}$ minimize $\mathcal L_{\rm pop}$ on $\mathcal U_N$; it exists by compactness. Comparing the empirical objective at $\widehat\beta$ and $\beta_N^{\rm opt}$ yields
\[
\begin{aligned}
 \mathcal L_{\rm pop}(\widehat\beta)-\mathcal L^*
 &\le \mathcal L_{\rm pop}(\beta_N^{\rm opt})-\mathcal L^*
       +2\varepsilon_{\rm num}+2\varepsilon_{\rm samp}+\rho_K\\
 &\le L_S^2d_N^2+2\varepsilon_{\rm num}+2\varepsilon_{\rm samp}+\rho_K.
\end{aligned}
\]
The second inequality uses the feasible best uniform approximation to $\beta^*$ and the centered-noise identity~\eqref{eq:sir_population_identity}. Thus the square on $d_N$ uses both realizability of the true control and conditional centering of the noise; it is not the generic Lipschitz objective estimate. Substituting the numerical and sampling bounds gives, with probability at least \(1-\alpha\),
\begin{equation}\label{eq:sir_explicit_bound}
 \begin{aligned}
 \mathcal L_{\rm pop}(\widehat\beta)-\mathcal L^*
 &\leq\underbrace{L_S^2d_N^2}_{\text{Approximation bound}}
      +\underbrace{4(1+\sigma)E_h}_{\text{Numerical bound}}\\
 &\quad+\underbrace{2\eta_{M,N}(r,\alpha)}_{\text{Generalization bound}}
      +\underbrace{\rho_K}_{\text{Optimization gap}}.
 \end{aligned}
\end{equation}
Here \(L_S\) is given by \eqref{eq:sir_sensitivity_constants}, \(d_N\) by
\eqref{eq:sir_approximation_rate}, \(E_h\) by \eqref{eq:sir_euler_bound},
and \(\eta_{M,N}\) by \eqref{eq:sir_uniform_sampling}. The estimate concerns
the time-averaged squared error of the infected trajectory through
\eqref{eq:sir_population_identity}. It does not, by itself, give a rate for
recovering the transmission coefficient.

\begin{remark}[Identifiability for positive initial populations]
\label{rmk:sir_identifiability}
Suppose additionally that $S_0>0$ and $I_0>0$. With the initial state and
$\gamma$ known, the map $\beta\mapsto I_\beta$ is injective on
$\mathcal B_{L_0}$. Under the stated realizability and conditional-centering
assumptions, $\beta^*$ is therefore the unique minimizer of
$\mathcal L_{\rm pop}$. Indeed, for every admissible $\beta$,
\[
 \begin{aligned}
 S_\beta(t)&=S_0\exp\!\left(-\int_0^t\beta(s)I_\beta(s)\,\mathrm ds\right)>0,\\
 I_\beta(t)&=I_0\exp\!\left(\int_0^t[\beta(s)S_\beta(s)-\gamma]\,\mathrm ds\right)>0.
 \end{aligned}
\]
The infected trajectory is continuously differentiable and determines
\[
 R_\beta(t)=R_0+\gamma\int_0^t I_\beta(s)\,\mathrm ds,
 \qquad S_\beta(t)=1-I_\beta(t)-R_\beta(t),
\]
and hence
\[
 \beta(t)=\frac{\dot I_\beta(t)+\gamma I_\beta(t)}
                    {S_\beta(t)I_\beta(t)}.
\]
Any population minimizer has the same infected trajectory as $\beta^*$
by \eqref{eq:sir_population_identity} and continuity, so injectivity gives
uniqueness. If $\widehat\beta_j\in\mathcal B_{L_0}$ and
$\mathcal L_{\rm pop}(\widehat\beta_j)-\mathcal L^*\to0$, compactness and
continuity imply that every subsequential limit is $\beta^*$. Thus
$\|\widehat\beta_j-\beta^*\|_\infty\to0$. This applies under the refinement
and vanishing global-gap hypotheses of Theorem~\ref{thm:gamma_conv},
on the stated sample event for empirical objectives. A quantitative
coefficient-recovery rate would additionally require an inverse
stability estimate.
\end{remark}

The bound separates approximation, numerical, sampling, and optimization
errors under the stated assumptions. Increasing \(N\) decreases the
interpolation error bound and changes the parameter-covering estimate.
Reducing \(h\) decreases the displayed state-error estimate, while
increasing \(M\) decreases the sampling estimate for a fixed class and
covering radius.
Running an optimizer longer is useful for this bound insofar as it reduces
the global gap \(\rho_K\). No architecture-independent width rule or universal
optimizer iteration rate follows. The exponential factor in \(L_S\) can be
large on long intervals, so the uniform bound may be conservative.

\subsection{Numerical Illustration}
\label{subsec:sir_numerical_illustration}

We use synthetic observations to examine these distinctions, with
\(T=12\), \(\gamma=0.2\), \(x_0=(0.98,0.02,0)\), and
\[
 \beta^*(t)=0.5+0.2\sin(2\pi t/12),\qquad
 0.2\leq\beta\leq0.8,\quad L_0=0.12.
\]
All training trajectories use Euler with linear reconstruction. Numerical
reference trajectories use DOP853 with relative tolerance
\(2\times10^{-12}\), absolute tolerance \(2\times10^{-14}\), and maximum
step \(0.15\), splitting integration at spline knots. Reference integrals
use composite Gauss--Legendre quadrature. These computations provide
reference values rather than tolerance-based error certificates.
All reported slopes are fitted on the four finest resolutions.
All random quantities use the recorded base seed (20261011), with
separate streams for the sampling, learning, and scalar experiments.
The supplied Python program reproduces the observations, fitted
coefficients, and figures. Only the infected component is observed; the
initial state and recovery rate remain known throughout.

The numerical implementation was reviewed with assistance from ChatGPT
(OpenAI). The reported numerical results and figures were computed from the
specified models and numerical procedures using the accompanying Python code.

\begin{figure}[tbp]
 \centering
 \includegraphics[width=\textwidth]{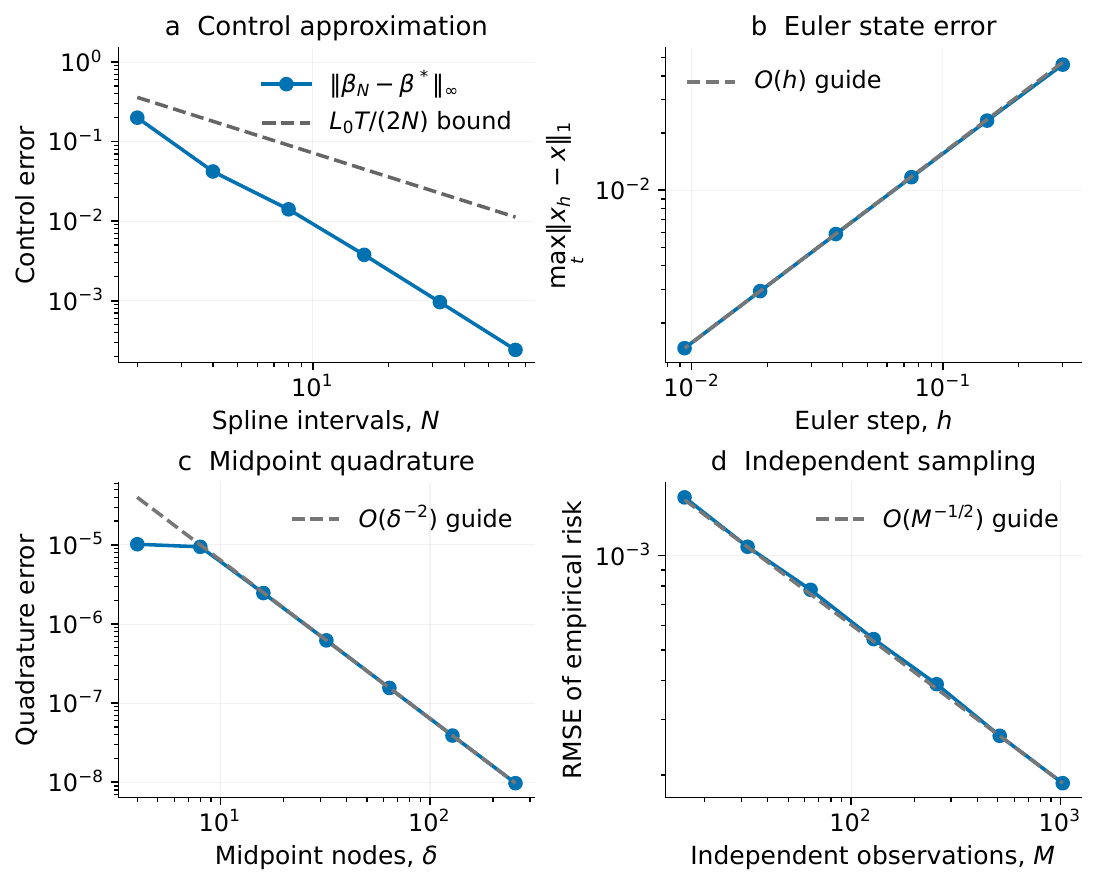}
 \caption{Separate error diagnostics for the SIR problem. (a) Control
 interpolation error and its Lipschitz bound. (b) Euler state error on a
 dense reference grid. (c) Midpoint integration error for a fixed loss.
 (d) Empirical-risk root mean squared error over independently generated
 datasets for a fixed control. Dashed lines in
 (b)--(d) indicate reference slopes, not theoretical upper bounds.}
 \label{fig:sir_error_diagnostics}
\end{figure}

Figure~\ref{fig:sir_error_diagnostics} shows separate approximation, Euler,
quadrature, and sampling diagnostics. Interpolating \(\beta^*\) with
\(N=2,4,8,16,32,64\) gives
orders \(1.96\) for the control error, \(1.99\) for the state error, and
\(3.99\) for the squared infected-trajectory loss. The smooth truth satisfies the additional regularity needed
for the higher approximation order. State-error maxima are computed on a
finite grid of \(12001\) equally spaced reference times. With the true
control fixed,
Euler steps from \(0.3\) to \(0.009375\) give state-error order \(1.00\).
For the fixed admissible comparison control
\(0.44+0.14\cos(2\pi t/12+0.35)\), midpoint integration of the noiseless
squared loss has observed order \(2.00\). Independently drawing uniform
observation times and noise \(\xi\sim\operatorname{Unif}[-0.003,0.003]\)
gives a distinct sampling experiment. Over \(2000\) independent datasets
per sample size, the empirical-risk root mean squared error scales as
\(M^{-0.51}\). Predictions in this fixed-control sampling diagnostic use
the reference trajectories. This diagnostic does not measure a uniform
generalization bound for trained controls.

\begin{figure}[tbp]
 \centering
 \includegraphics[width=\textwidth]{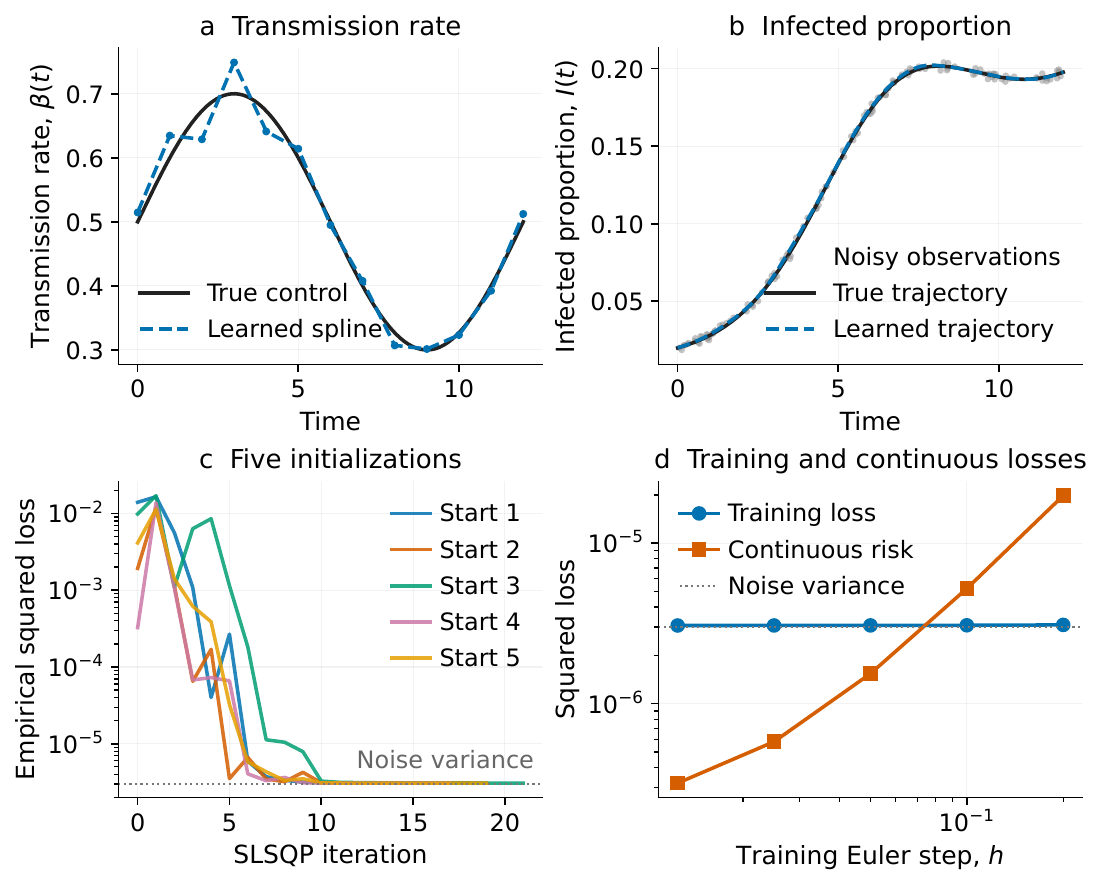}
 \caption{Learning from infected-component observations. (a) True and fitted
 transmission rates. (b) Noisy observations and continuous reference
 trajectories. (c) Empirical losses from five feasible initializations;
 the noise variance is a reference level, not a finite-sample lower bound.
 (d) Discrete training loss and continuous noiseless risk after refitting
 at different Euler steps. Panels (a)--(c) use \(N=12\), \(M=128\), and
 \(h=0.025\).}
 \label{fig:sir_learning}
\end{figure}

For learning, we use \(M=128\) noisy observations, \(N=12\) spline
intervals, and \(h=0.025\). Sequential least-squares programming (SLSQP)
minimizes the empirical loss under the
nodal amplitude and difference constraints, using derivatives of the
discrete Euler recursion. Five feasible initializations reach losses of
approximately \(3.07108\times10^{-6}\) in \(17\)--\(21\) iterations.
The selected fit has continuous root mean squared errors
\(7.61\times10^{-4}\) for \(I\) and \(1.60\times10^{-2}\) for \(\beta\).
Agreement among these initializations describes this instance and does not establish
global optimality. Figure~\ref{fig:sir_learning} shows the fitted
coefficient, observations, and optimization histories.

Refitting the same data and spline class at different Euler steps exposes
a discrepancy that the training loss alone does not reveal
(Table~\ref{tab:sir_solver_refinement}). Reducing \(h\) from \(0.2\) to
\(0.0125\) changes the empirical loss by less than one percent, while
reducing the continuous noiseless risk by a factor of approximately \(62\).
In this synthetic example, the fitted coefficient partly compensates for
the solver error. An accurate discrete fit can therefore conceal a larger
error when the learned coefficient is evaluated in the continuous model.

\begin{table}[tbp]
 \centering
 \caption{Solver refinement with the same observations and spline class.
 Continuous risk is the time-averaged squared error from the noiseless
 infected trajectory, evaluated by reference integration.}
 \label{tab:sir_solver_refinement}
 \begin{tabular}{ccc}
 \hline
 Euler step & Training loss & Continuous risk\\
 \hline
 \(0.2\)    & \(3.09565\times10^{-6}\) & \(1.98020\times10^{-5}\)\\
 \(0.1\)    & \(3.07869\times10^{-6}\) & \(5.21820\times10^{-6}\)\\
 \(0.05\)   & \(3.07344\times10^{-6}\) & \(1.53805\times10^{-6}\)\\
 \(0.025\)  & \(3.07108\times10^{-6}\) & \(5.79238\times10^{-7}\)\\
 \(0.0125\) & \(3.06967\times10^{-6}\) & \(3.19505\times10^{-7}\)\\
 \hline
 \end{tabular}
\end{table}

Finally, a scalar instance allows us to bound the optimization gap independently of the optimizer.
We take \(T=2\), \(x_0=(0.9,0.1,0)\), \(\gamma=0.2\), and constant
\(b\in[0.2,0.8]\), with true value \(b^*=0.5\). There are \(64\)
independent uniform observation times, uniform noise on \([-0.01,0.01]\),
and \(256\) Euler steps, so \(h=1/128\). For every admissible \(b\), the
Euler update gives
\[
 I_j\leq0.1(1+0.6h)^j\leq0.1e^{1.2},\qquad 0\leq j\leq256.
\]
The exact infected trajectory has the same envelope, and linear
reconstruction preserves the numerical nodal bound. Comparing two
constant controls in the Euler recursion gives the state-map Lipschitz
constant
\[
 C_{\rm state}=\frac12\frac{(1+1.6h)^{256}-1}{1.6}.
\]
Writing $B_I=0.1e^{1.2}$, predictions lie in $[0,B_I]$ and observations lie in $[-0.01,B_I+0.01]$. Hence $|2(p-y)|\le2(B_I+0.01)$ throughout this range, which gives the squared-loss Lipschitz constant \(2(0.1e^{1.2}+0.01)\) in the infected component. Multiplying by
\(C_{\rm state}\) therefore bounds the Lipschitz constant of the discrete
empirical mean squared loss \(F(b):=\mathcal L_{64,1/128}(b)\) by
\[
 C_F=(0.1e^{1.2}+0.01)
      \frac{(1+1.6h)^{256}-1}{1.6}\approx4.92726.
\]
On a uniform parameter grid of spacing \(\Delta b\),
\begin{equation}\label{eq:sir_scalar_grid_certificate}
 \inf_{b\in[0.2,0.8]}F(b)
 \geq\max\{0,\min_{b\ \mathrm{on\ grid}}F(b)-C_F\Delta b/2\}.
\end{equation}
The grid includes both endpoints, so every admissible parameter lies
within \(\Delta b/2\) of a node. Lipschitz continuity therefore controls
the objective throughout the interval, independently of the optimizer.
With \(2^{17}\) grid intervals, bounded scalar minimization returns
\(\widehat b=0.49429818\). Its loss is \(3.09073608\times10^{-5}\),
and the grid lower bound is \(1.96297558\times10^{-5}\). Thus the scalar
gap \(\rho_{\rm sc}:=F(\widehat b)-\min_{b\in[0.2,0.8]}F(b)\) has a
bound of approximately \(1.12776\times10^{-5}\).

The true constant is admissible, so the approximation term is zero.
For \(\mathcal C=\{\widehat b,b^*\}\), the numerical and sampling
discrepancies used in the comparison are
\begin{align*}
 A_{\rm num}&=\sum_{b\in\mathcal C}|\mathcal L_{\rm pop}(b)-\mathcal L_h(b)|,\\
 A_{\rm samp}&=\sum_{b\in\mathcal C}|\mathcal L_h(b)-F(b)|.
\end{align*}
Inserting the intermediate objectives at these two controls gives the
population-excess estimate \(A_{\rm num}+A_{\rm samp}+\rho_{\rm sc}\).
Here \(A_{\rm num}=2.21935\times10^{-8}\) and
\(A_{\rm samp}=4.79978\times10^{-6}\). Together with the gap bound they
give a comparator budget \(1.60996\times10^{-5}\), exceeding the observed
population excess \(5.77604\times10^{-7}\). These two discrepancies are
reference evaluations at the comparison controls, not inferred uniform
bounds. Formula~\eqref{eq:sir_scalar_grid_certificate} is analytic for the
discrete objective; its floating-point evaluation and the continuous
reference integrals are not interval-arithmetic certificates. No comparable
global certificate is claimed for the preceding multivariate SLSQP fits.

\section{Conclusion and Discussion}
We have analyzed time-dependent learned terms coupled to prescribed mechanistic ODEs, a specialization of universal differential equations \cite{rackauckas2020universal,yin2023optimal}. The objective error bound separates neural approximation, numerical integration, statistical estimation, and optimization errors within a fixed compact control class. Deterministic quadrature provides an alternative form of loss evaluation with its own estimate.

Constrained nodal ReLU realizations and Euler integration give explicit approximation and state-error estimates. The statistical bound accounts for the complete control-to-trajectory map, and the variational analysis retains the changing feasible domains. Under joint refinement and vanishing global objective gaps, minimum values converge and approximate minimizers approach the continuous optimal set. The common regularity bound is central to this compactness argument.

The SIR example makes the constants explicit on the invariant population simplex. For positive initial susceptible and infected populations, a known initial state, and a known recovery rate, the population minimizer is identifiable; vanishing excess risk then implies uniform convergence of the transmission rate within the admissible class. The numerical experiments isolate interpolation, integration, quadrature, and fixed-control sampling discrepancies. Refitting the same observations at several solver resolutions shows that nearly unchanged training loss can accompany a much larger continuous-trajectory error on a coarse grid. The scalar experiment further connects a lower bound for the discrete minimum with an upper bound for the optimization gap.

The global optimization gap remains an explicit condition for the general nonconvex problem. Quantitative recovery of coefficients requires an inverse stability estimate, while identification in other models depends on their observation maps. Higher-order solvers require regularity compatible with the learned control. Further extensions include adaptive refinement, broader control classes, dependent observations, and state-dependent learned terms beyond the stability estimates proved here. PDE-constrained problems additionally require suitable state-space compactness and uniform discretization estimates.

\section*{Acknowledgements}
The authors would like to thank the referees for their careful reading and constructive comments.

\section*{Funding}
Research partially supported by the Guangdong Provincial Key Laboratory of Mathematical and Neural Dynamical Systems (2024B1212010004), the National Natural Science Foundation of China (12571544, 12220101001), the National Key R\&D Program of China (2022YFA1003704), the China Postdoctoral Science Foundation (No.\ 2023T160504), the Shaanxi Institute of Fundamental Science Project (No.\ 23JSQ028, 25JSY034), and the Xi'an Jiaotong University Basic Research Operating Expenses (No.\ xyz012023045).

\section*{Data and code availability}
All observations used in the numerical illustration are synthetic and are generated from the SIR systems and noise distributions specified in Section~\ref{sec:sir_example}. The numerical source code, fixed random seed, generated data, and figure-generation routines are provided in the ancillary files accompanying this preprint. No clinical or personal data were used.

\section*{Declaration of generative AI and AI-assisted technologies}
ChatGPT (OpenAI) was used to assist with language editing and review of the
mathematical exposition and numerical implementation. The authors are
responsible for all mathematical arguments, numerical results, and conclusions
presented in this article.

\bibliographystyle{elsarticle-num}
\bibliography{references}
\end{document}